\documentclass[onefignum,onetabnum]{siamart171218}

\usepackage{lipsum}
\usepackage{amsfonts}
\usepackage{graphicx}
\usepackage{epstopdf}
\usepackage{algorithmic}
\ifpdf
  \DeclareGraphicsExtensions{.eps,.pdf,.png,.jpg}
\else
  \DeclareGraphicsExtensions{.eps}
\fi

\newsiamremark{remark}{Remark}
\newsiamremark{hypothesis}{Hypothesis}
\crefname{hypothesis}{Hypothesis}{Hypotheses}
\newsiamthm{claim}{Claim}

\headers{Discrete Stream(line) Method}{Kiril S. Shterev}

\title{On convective terms approximation approach that corresponds to pure convection\thanks{
\funding{This work was funded by the Bulgarian NSF under Grant DN-02/7-2016.}}}

\author{Kiril S. Shterev\thanks{Institute of Mechanics, Bulgarian Academy of Sciences, Acad. G. Bonchev Str., Block 4, Sofia 1113, Bulgaria
  (\email{kshterev@imbm.bas.bg}, \url{www.imbm.bas.bg}).}}

\usepackage{amsopn}

\ifpdf
\hypersetup{
  pdftitle={On convective terms approximation approach that corresponds to pure convection},
  pdfauthor={Kiril S. Shterev}
}
\fi

\usepackage{lineno}

\usepackage{graphicx}
\usepackage{caption}
\usepackage{subcaption}

\usepackage{amssymb}

\usepackage{amsmath}

\usepackage{geometry}

\usepackage{multirow}
\usepackage{ctable}

\usepackage{placeins}

\usepackage[utf8]{inputenc}
\usepackage[T1, T2A]{fontenc}

\begin{document}

\maketitle

\begin{abstract}
Recent decades are put lots of efforts to develop a higher-order scheme for convective terms approximation that is stable and reliable. The idea presented here is that approximation approach has to correspond to the physical phenomenon described by approximated terms. Pure convection (advection) that is described by convective terms is transporting a property along the streamline, and the information propagation is unidirectional, i.e., transported property depends on previous values along the streamline but does not depend on the next ones. The proposed approach represents streamlines on mesh as discrete streamlines and is called Discrete Stream(line) Method (DStreaM). A discrete streamline here is represented as a narrow triangle with one vertex of the approximated node and two others neighbor upstream nodes. Discrete streamlines are orientated using local flow direction as skew upwind schemes. DStreaM corresponds to pure convection. Here are considered standard test problems: advection of a step profile, advection of a double-step profile, advection of a sinusoidal profile, and Smith and Hutton problem. DStreaM solutions were compared with upwind-first order scheme and second-order Total Variation Diminishing (TVD) schemes with limiters Min-Mod, QUICK, and SUPERBEE solutions. DStreaM demonstrated second-order accuracy and rapid convergence. Upwind and DStreaM need 2 or 4 iterations to reach a final solution while TVD schemes need from 15 to 93.5 more iterations. DStreaM approach looks promising for calculation of convective-dominated problems because it approximates naturally first derivatives and is straightforwardly applicable as a meshfree method or on unstructured meshes.
\end{abstract}



\section{Introduction}
Development of higher-order convective terms approximation scheme is an important problem of present-day science. Navier-Stokes (NS) equations are one of the dominant equations used in Computational Fluid Dynamic (CFD) to model fluid flows. NS system of equations describes unsteady, viscous, incompressible/compressible, and heat-conductive fluid. They contain unsteady, convective, and diffusion terms. As the diffusion terms could be easily approximated with at least second-order spatial accuracy, convective terms can be approximated unconditionally with upwind first-order scheme. Upwind first-order scheme obtains correct results for all speeds and parameters, but it is computationally expensive. The corresponding second-order scheme that obtains correct physical results for all speeds and parameters as upwind first-order scheme does not exist. Up to now are followed two routes to develop higher-order convective terms approximation scheme: One way is to use higher-order schemes that are basically derived as one-dimensional scheme and applied for multidimensional cases and help in reducing numerical diffusion. Total Variation Diminishing (TVD) schemes are a typical representative that follow this route. They are designed to prevent the undesirable oscillatory behavior of higher-order schemes. In TVD schemes, the tendency towards oscillation is counteracted by adding an artificial diffusion fragment or by adding a weighting towards upstream contribution. In the literature early schemes, based on these ideas were called flux corrected transport (FCT) schemes: see Boris and Book \cite{Boris1973}, \cite{Boris1976}. Further, works by van Leer \cite{vanLeer1974}, \cite{VanLeer1977_III}, \cite{VanLeer1977_IV}, \cite{vanLeer1979}, Harten \cite{Harten1983}, \cite{Harten1984}, Sweby \cite{Sweby1984}, Roe \cite{Roe1985}, Osher, and Chakravarthy \cite{Osher1984}, Zijlema \cite{Zijlema1996}, Arora and Roe \cite{Arora1997}, \v{C}ada and Torrilhon \cite{Cada2009}, Ferreira et al. \cite{Ferreira2012}, Zhang et al. \cite{Zhang2015}, Kriel \cite{Kriel2017}, and many others have contributed to the development of present-day TVD schemes. Another route is to use a skew upwind scheme that reduces numerical diffusion because of their multidimensional nature. They take into account flow direction (for 2D and 3D flows) and determine neighbor points according to it, in some sense, they follow the streamlines. First Raithby presented skew upwind schemes SUDS and SUWDS in \cite{Raithby1976_1} and \cite{Raithby1976_2}. These schemes tend to involve the neighboring diagonal cells as well as those opposite to each cell face. Also, within the finite volume context, Hassan, Rice, and Kim \cite{Hassan1983} proposed a mass-flow-weighted skew upwind scheme as an improvement over the conventional skew-upwind scheme. The scheme ensures a reduction of numerical instability and numerical diffusion errors. In the finite element context, Schneider and Raw \cite{Schneider1986} proposed an upwind procedure that accounts for the directionality of the flow field through a skewed approach, while simultaneously precluding the possibility of negative coefficients. Patel, Cross, and Markatos \cite{Patel_et_all_1988} retained the general objectives of Raithby scheme and proposed CUPID scheme. CUPID scheme accounts for the flow angles at the corners of a cell rather than at the cell faces, as does Raithby's SKEW schemes. Skew upwind weighted differencing (SUWD) scheme was proposed by Busnaina, Zheng, and Sharif \cite{Busnaina1991}. Here the convected quantity at a cell face is approximated by a weighted average of two quantities obtained by using linear interpolation at two locations upstream of the resultant velocity vector at the cell face. Carey, Scanlon, and Fraser \cite{Carey_and_Scanlon_and_Fraser_1992}, \cite{Carey_and_Scanlon_and_Fraser_1993} created SUCCA scheme that concentrates attention at the corners of finite-volume cells. Darwish and Moukalled \cite{Darwish_and_Moukalled_1996} presented STOIC scheme (Second- and Third-Order Interpolation for Convection) that is a high-resolution scheme developed and implemented in the context of the normalized variable formulation methodology developed by Leonard \cite{Leonard1988}. Ogedengbe and Naterer presented convective upwind scheme called NISUS (Non-Inverted Skew Upwind Scheme) in \cite{Ogedengbe2004} and \cite{Ogedengbe2006}. The main benefit of NISUS lies in avoiding the costly inversion of the upwind coefficient matrix, without any significant loss of accuracy. Karadimou and Markatos retain the objectives of CUPID and SUCCA schemes and reformulate the convection terms in the momentum and scalar conservation equations in a way to treat 3D flows. They presented SUPER (Skew Upwind and Corner Algorithm) scheme in \cite{Karadimou2012}. After all, both routes obtain more accurate results than first-order upwind scheme, but they can give rise of oscillations or under/over-shoots, especially in regions of strong gradients when are applied to calculate fast compressible flows. TVD schemes particularly obtain excellent results for the steady state, slow and moderate flows. Slow and moderate flows are unsteady but close to steady one for lower Reynolds numbers. When unsteady fluid flows is fast TVD schemes are not time accurate, see Chung \cite{Chung_2002}. Furthermore, if Reynolds number is higher, even when fluid flow is steady, TVD schemes can obtain an unphysical oscillated solution. In such cases, the solution can be obtained using the first-order upwind scheme using significant computational resources. Van Leer in his review article \cite{van_Leer_2006} writes that American aeronautical community is skeptical about the cutting-edge algorithm, and they prefer the advent of high-performance computing and promise of massively parallel computing instead of the development of higher order approximation schemes. The need for a new approach is evident.\\\indent
The approximation scheme should correspond to process that approximated term(s) describe. Diffusion terms describe exchange of concentration of specific property from one place to another. Diffusion depends on concentration and does not depend on flow direction, i.e., does not depend on velocity. Information is exchanged in all directions according to a specific node. One of the simplest approximation schemes for diffusion terms is the central difference, it is unconditionally stable when approximate diffusion terms, it is a second-order scheme, and it takes information from neighbor nodes in all directions. The convection opposite to diffusion depends on flow direction and does not depend on concentration. Convection terms describe transportation of property ($\phi$) along the streamline. They contain information about streamline (velocities in convective terms) and previous values of $\phi$ (see equation (\ref{eq1})) along the streamline, but they do not contain information for $\phi$ further along the streamline. Approximation scheme that corresponds to convection has to contain upwind information along the streamline that passes through the approximated node and no downstream information about $\phi$. TVD schemes use downstream information about $\phi$ that corresponds to diffusion, which is not included as information in convective terms. This is their problem and the reason for their limited application.\\\indent
In this paper is presented approximation approach for convective terms that use only upstream points and orients approximation scheme along streamlines that makes it of a group of skew upwind schemes. Four standard two-dimensional pure convection (advection) steady test cases were used to compare the relative performance of proposed Discrete Stream(line) Method (DStreaM) with standard schemes: upwind first-order and second-order TVD schemes. TVD schemes are well studied and widely used in Computational Fluid Dynamic (CFD) codes that make them appropriate schemes for presented comparisons. Here are used TVD scheme with limiters Min-Mod \cite{Roe1985}, QUICK \cite{Leonard1988} and SUPERBEE \cite{Roe1985}. Sweby \cite{Sweby1984} has given necessary and sufficient conditions for a scheme to be TVD second-order. The Min-Mod limiter function exactly traces the lower limit of the TVD region, SUPERBEE limiter follows the upper limit, and QUICK limiter is between them. The presented test problems are advection of a step profile, advection of a double-step profile, advection of a sinusoidal profile, and Smith and Hutton problem. Obtained results, show the second-order spatial accuracy of DStreaM and significant reduction of iterations from 15 to 93.5 times according to TVD schemes.
\section{Considerations on convective and diffusion terms information propagation}
\begin{figure}[htb!]
	\centering
    \begin{subfigure}[b]{0.48\textwidth}
		\centering
        \includegraphics[width=\textwidth]{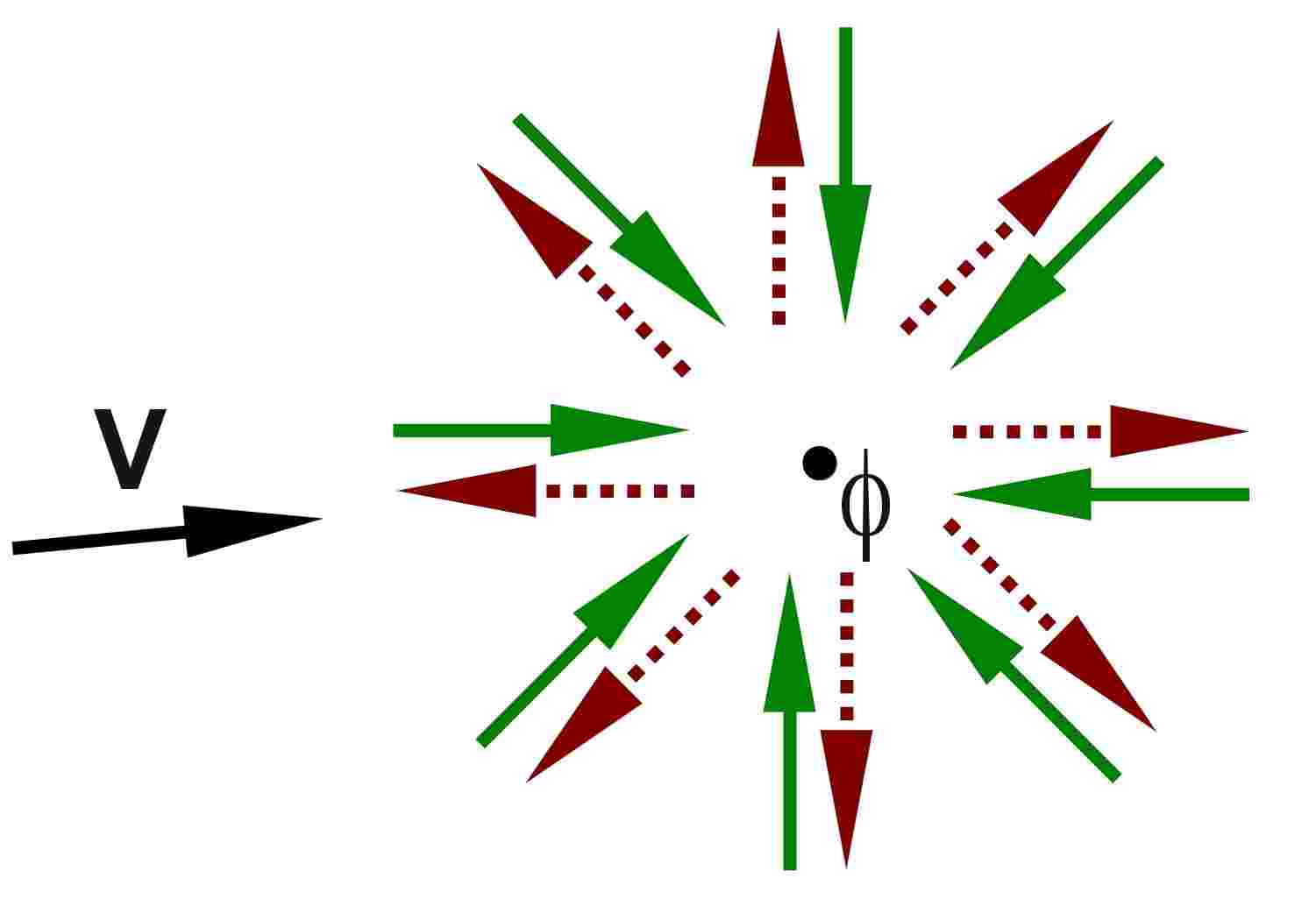}
        \caption{}
        \label{Propagation_information:diffusion_terms}
    \end{subfigure}
    \begin{subfigure}[b]{0.48\textwidth}
		\centering
        \includegraphics[width=\textwidth]{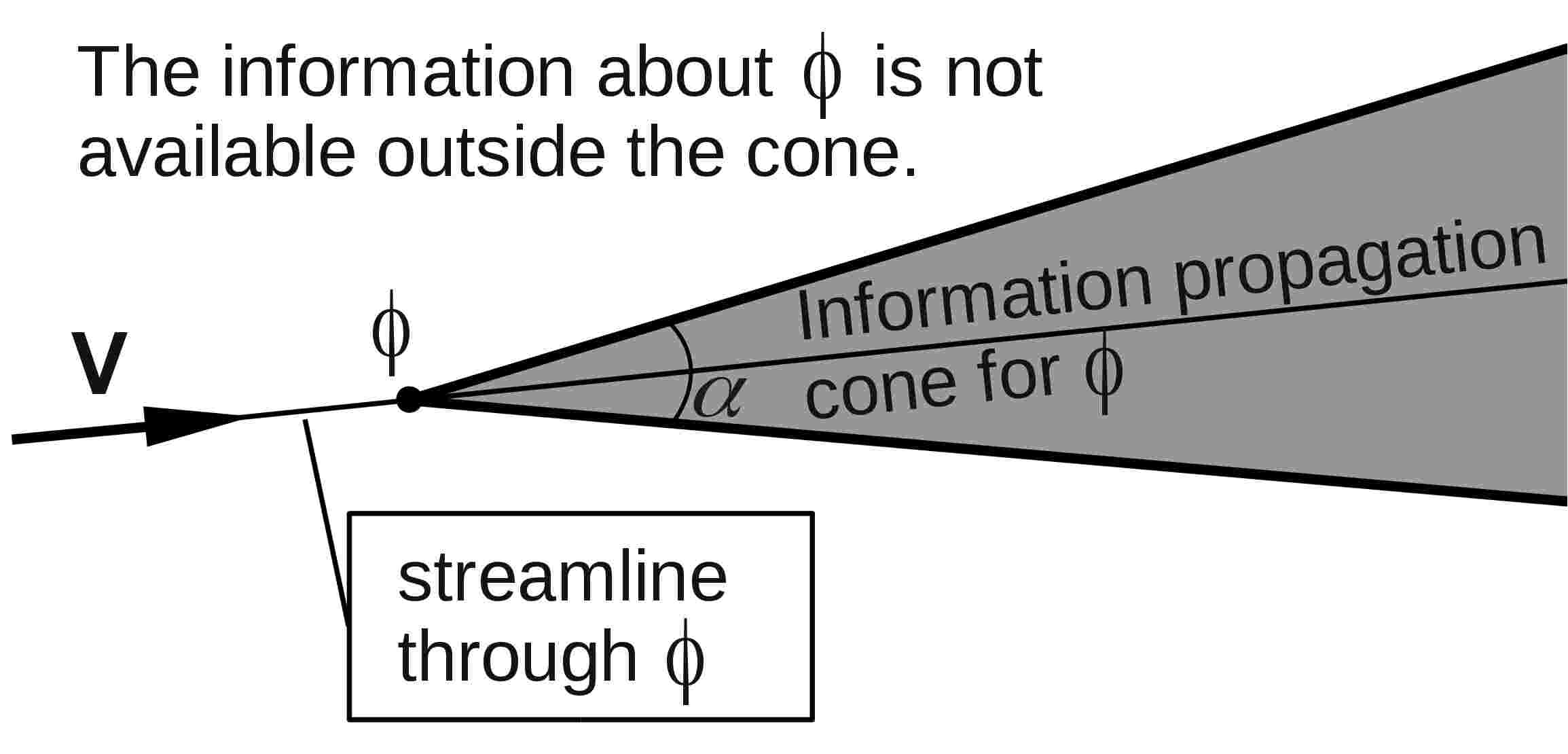}
        \caption{}
        \label{Propagation_information:convection_terms}
    \end{subfigure}
    \caption{Schematic information propagation of diffusion dominated problem (a) and convection dominated problem (b).}
    \label{Propagation_information}
\end{figure}

\begin{figure}[htb!]
	\centering
    \includegraphics[width=0.5\textwidth]{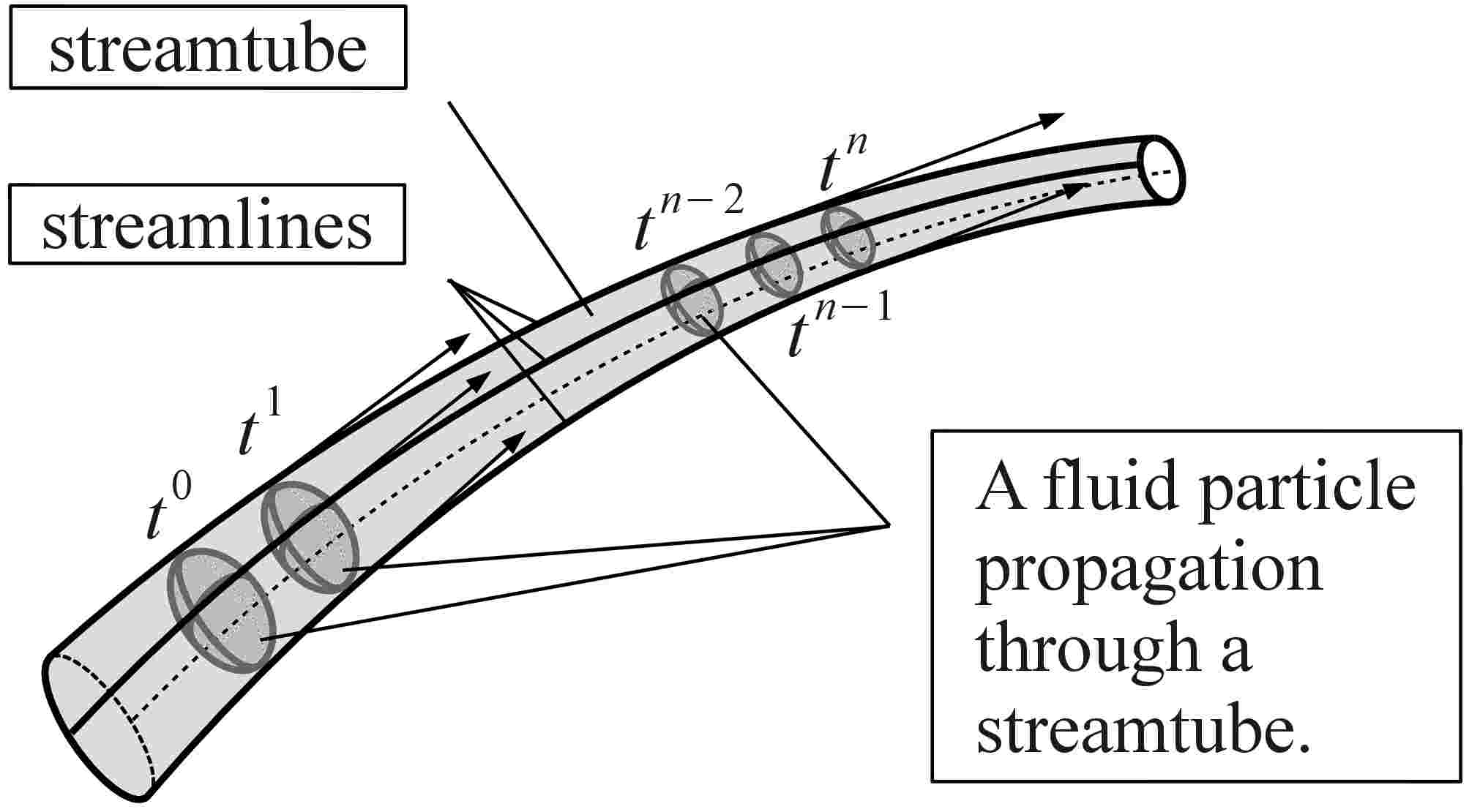}
    \caption{Transport of fluid particle through streamtube.}
    \label{Transport_of_fluid_particle_through_streamtube}
\end{figure}

\begin{figure}[htb!]
	\centering
    \includegraphics[width=0.8\textwidth]{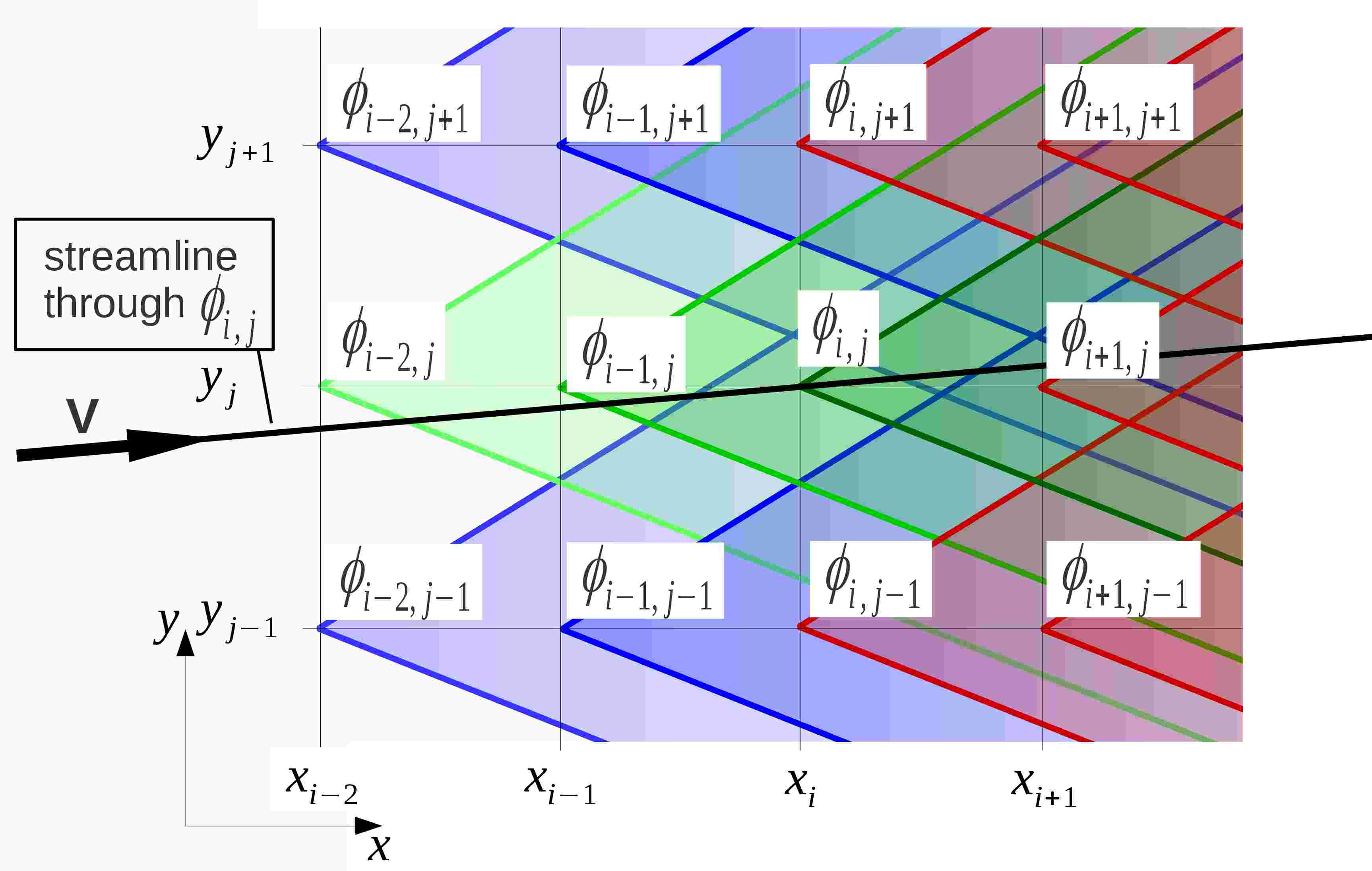}
    \caption{Information propagation cones of convection dominated problem at mesh nodes.}
    \label{Propagation_information_cones}
\end{figure}

Here are presented more detailed considerations about the pure convection and diffusion processes that are important for further explanations of the proposed approximation approach.\\\indent
The steady convection-diffusion equation for a general property $\phi$ is:
\begin{equation}
	\underbrace{\rho\frac{\partial(u \phi)}{\partial x} + \rho\frac{\partial(v \phi)}{\partial y}}_{\text{Convective terms}} = \underbrace{\Gamma \left( \frac{\partial^2\phi}{\partial x^2} + \frac{\partial^2\phi}{\partial y^2} \right)}_{\text{Diffusion terms}},
	\label{eq1}
\end{equation}
where $u$ and $v$ are velocities along x- and y-axis, $\rho$ is density, and $\Gamma$ is diffusion coefficient.\\\indent
The diffusion terms (see equation (\ref{eq1})) describe pure diffusion. Diffusion is the movement of molecules and/or atoms from a region of high concentration (or high chemical potential) to a region of low concentration (or low chemical potential). The information propagates in all directions independently of the velocity field. Fig. \ref{Propagation_information} (a) shows a diagram of directions that information propagates. Green continuous arrows denote the information that propagates to $\phi$ while red dashed arrows denote the information that propagates outward $\phi$. The information propagates in both directions. Therefore, the diffusion terms approximation scheme should include neighbor nodes in all directions. The central difference numerical scheme is second-order, matches perfectly to the diffusion terms properties, and approximates diffusion terms without need of any additional criteria to ensure convergence or prevent a non-physical solution.\\\indent
Fig. \ref{Propagation_information} (b) shows the information propagation of convective-dominated fluid flow schematically. The information about $\phi$ is available inside information propagation cone, but it is not available outside the cone including opposite or normal to velocity directions.\\\indent
The convective terms (see equation (\ref{eq1})) describe pure convection (advection). Advection can be considered as a motion of fluid particles through a streamtubes without exchange of properties between neighbor streamtubes. The transported fluid property is $\rho \phi$, and velocities in convective terms can be interpreted as fluid particle's streamline/path geometry, see Navier-Stokes system of equations derivation and considerations \cite{Loitsiansky2003}, \cite{White2005}. Fig. \ref{Transport_of_fluid_particle_through_streamtube} shows the motion of a fluid particle from time $t^{0}$ to time $t^{n}$ through a streamtube. The velocity field is known in all directions according to position of the fluid particle while the information relates transported property at time $t^{n}$ exist only in previous times of its travel up to $t^{0}$ and does not exist for next times as $t^{n+1}$ because they do not happen yet, and this kind of information is not included in convective terms. The unidirectional information propagation of transported property is the main characteristic of convective terms. When it is considered pure convection case, transported property of fluid particle at time $t^{n}$ is the same as at time $t^{n-1}$. Therefore, they are the same as transported property at time $t^{0}$: $\rho^{n} \phi^{n}=\rho^{n-1} \phi^{n-1}= ... =\rho^{0} \phi^{0}$. When the density is constant, it follows $\phi^{n}=\phi^{n-1}= ... =\phi^{0}$. That makes possible to determine an exact solution of a steady pure convection problem at first glance. The basic idea of the proposed approach uses convective terms unidirectional information transportation characteristic.
The propagation of information of convective-dominated problems in mesh is important to develop approximation that corresponds to modeled flow. Fig. \ref{Propagation_information_cones} shows propagation information cones of convective terms applied in nodes of a Cartesian uniform mesh. In ideal pure convection case each cone's angle $\alpha$ (see Fig. \ref{Propagation_information}, (b)) has to approach 0 ($\alpha \rightarrow 0$). Thus, will turn each cone into a curve and the node's value will be equal to value at the other side of the curve. That can be another node if the curve goes through it or the computational domain boundary. $\alpha \rightarrow 0$ corresponds to numerical scheme without false diffusion. False diffusion is information propagation normal to the streamline. Here cones represent information propagation of numerical solution of pure convection problems where false diffusion exists ($\alpha > 0$). The cones, in Fig. \ref{Propagation_information_cones}, can be separated into three groups according to the information that reaches $\phi_{i,j}$. The information from the first group reaches $\phi_{i,j}$. In considered case, these are $\phi_{i-1,j}$ and $\phi_{i-2,j}$. The information from the second group could reach $\phi_{i,j}$. It depends on the cone's angle ($\alpha$). The nodes of this group are $\phi_{i-2,j-1}$, $\phi_{i-2,j+1}$, $\phi_{i-1,j-1}$, and $\phi_{i-1,j+1}$. In present case information from node $\phi_{i-2,j-1}$ reaches $\phi_{i,j}$ while information from nodes $\phi_{i-2,j+1}$, $\phi_{i-1,j-1}$, and $\phi_{i-1,j+1}$ does not. And the information from the last group where are the nodes $\phi_{i,j-1}$, $\phi_{i,j+1}$, $\phi_{i+1,j-1}$, $\phi_{i+1,j}$, and $\phi_{i+1,j+1}$ does not reach $\phi_{i,j}$. The information propagation cone of $\phi_{i,j}$ rotated on 180 degrees is reversed information propagation cone and includes the nodes that information propagates to $\phi_{i,j}$, see Fig. \ref{Propagation_information_cones_that_reach_node_ij_UI}. According to these considerations, the approximation scheme of convective terms for the node $\phi_{i,j}$ can include only nodes that information reaches it and could not include nodes that information is impossible to achieve it.\\\indent
\section{Discrete Stream(line) Method (DStreaM)}
\begin{figure}[htb!]
	\centering
    \includegraphics[width=0.8\textwidth]{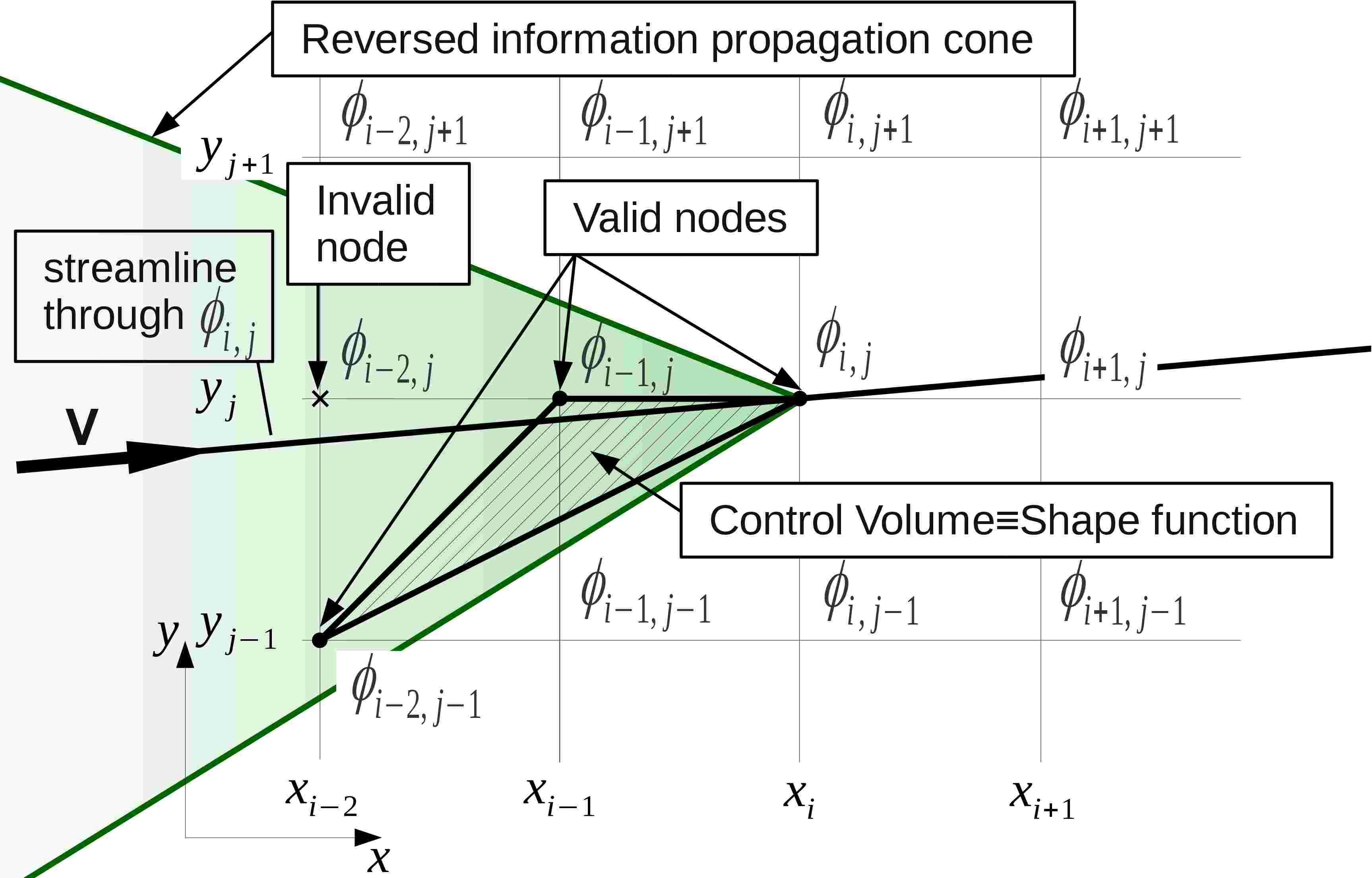}
    \caption{Reversed information propagation cone and shape function that corresponds to convective terms, and it is a discrete representation of a streamline through point $\phi_{i,j}$.}
    \label{Propagation_information_cones_that_reach_node_ij_UI}
\end{figure}

\begin{figure}[htb!]
	\centering
    \begin{subfigure}[b]{0.45\textwidth}
		\centering
        \includegraphics[width=\textwidth]{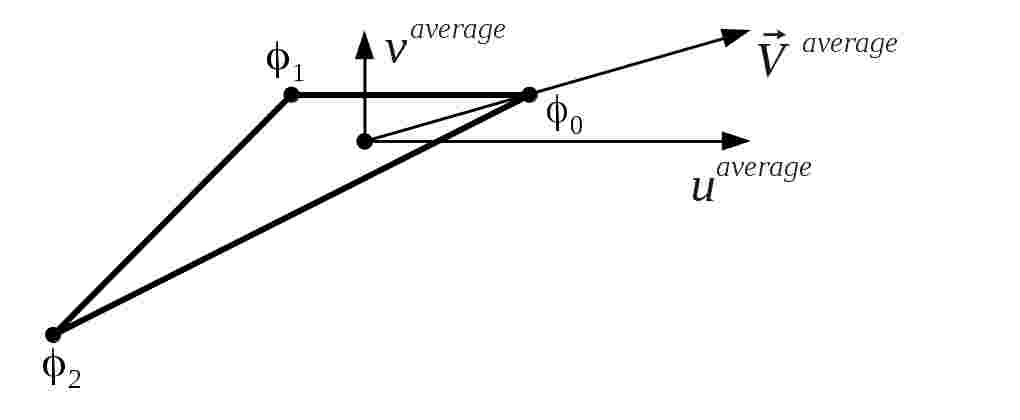}
        \caption{Shape function with average velocities.}
        \label{Shape_function}
    \end{subfigure}
    \begin{subfigure}[b]{0.45\textwidth}
		\centering
        \includegraphics[width=\textwidth]{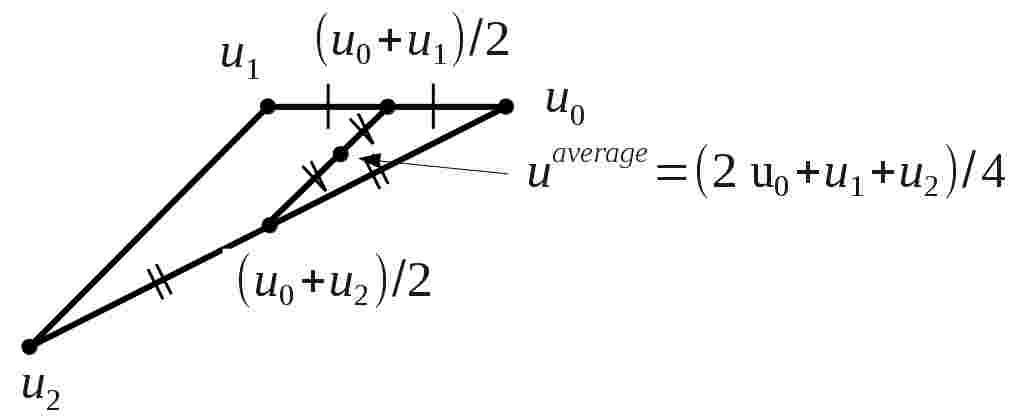}
        \caption{Horizontal average velocity approximation ($u^{average}$).}
        \label{Shape_function_u_average}
    \end{subfigure}
    \begin{subfigure}[b]{0.35\textwidth}
		\centering
        \includegraphics[width=\textwidth]{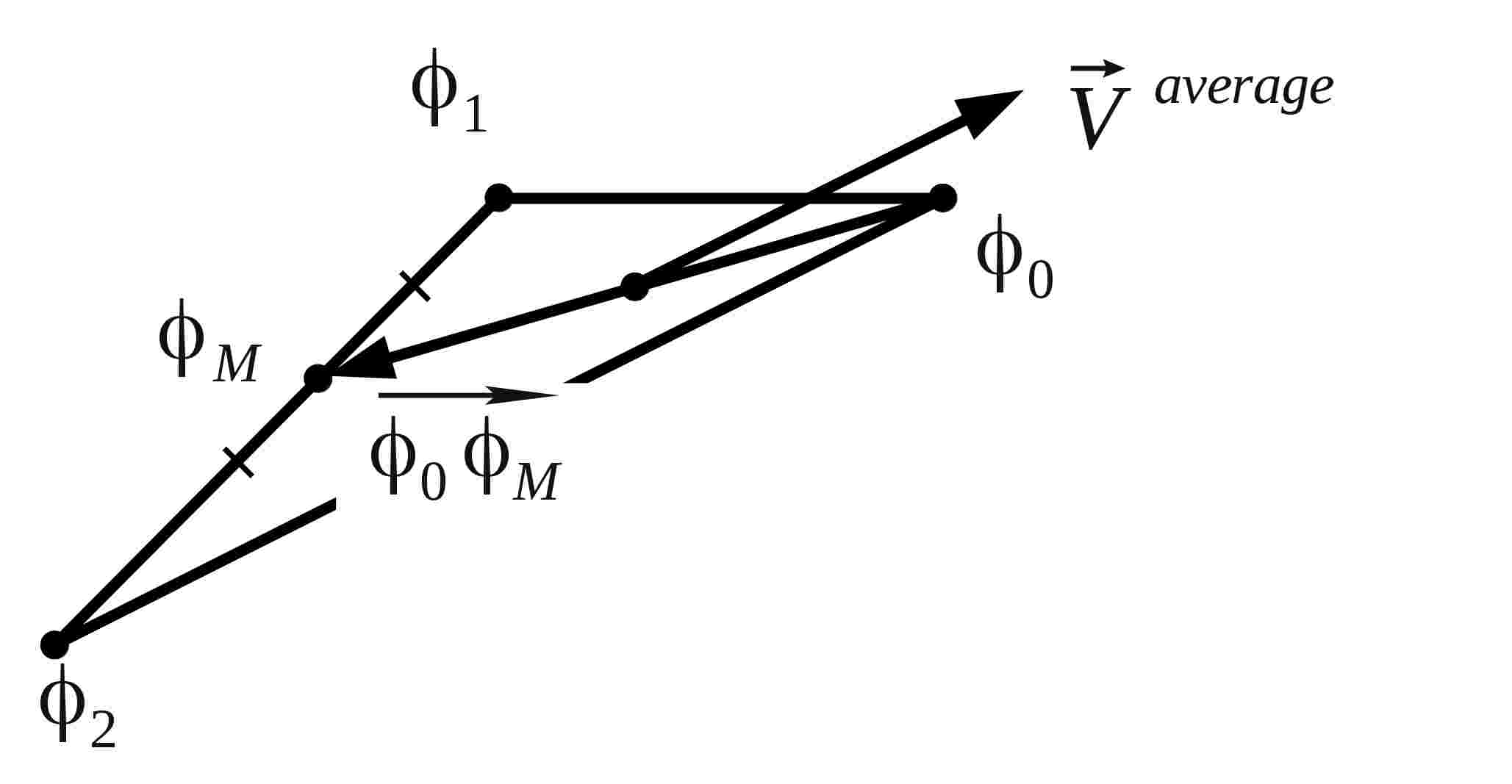}
        \caption{Vectors $\overrightarrow{V}^{average}$ and $\overrightarrow{\phi_0\phi_M}$ that are used to determine sector for approximation scheme.}
        \label{Shape_function_V_average_and_median}
    \end{subfigure}
    \caption{Shape function with average velocities (a), average velocities idea (b), and vectors that determine sector for approximation scheme (c).}
    \label{Shape_function_and_u_average}
\end{figure}

\begin{figure}
	\centering
	\begin{minipage}{0.32\textwidth}
		\ \\[2mm]
		\begin{subfigure}[b]{1.0\textwidth}
			\includegraphics[width=1\textwidth]{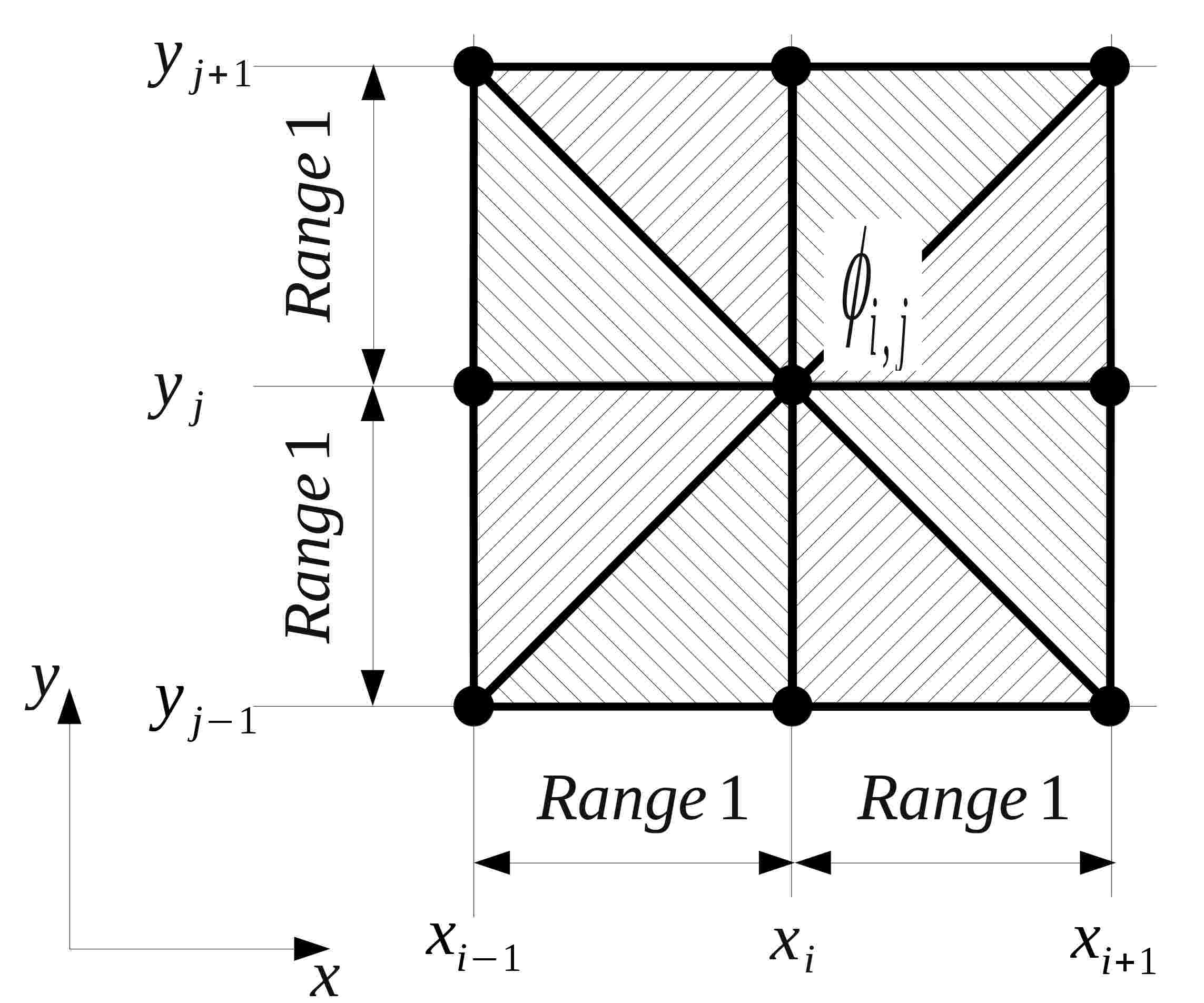}
			\caption{Shape functions for Range=1.}
		\end{subfigure}
	\end{minipage}
	\begin{minipage}{0.58\textwidth}
		\begin{subfigure}[b]{1.0\textwidth}
			\includegraphics[width=1\textwidth]{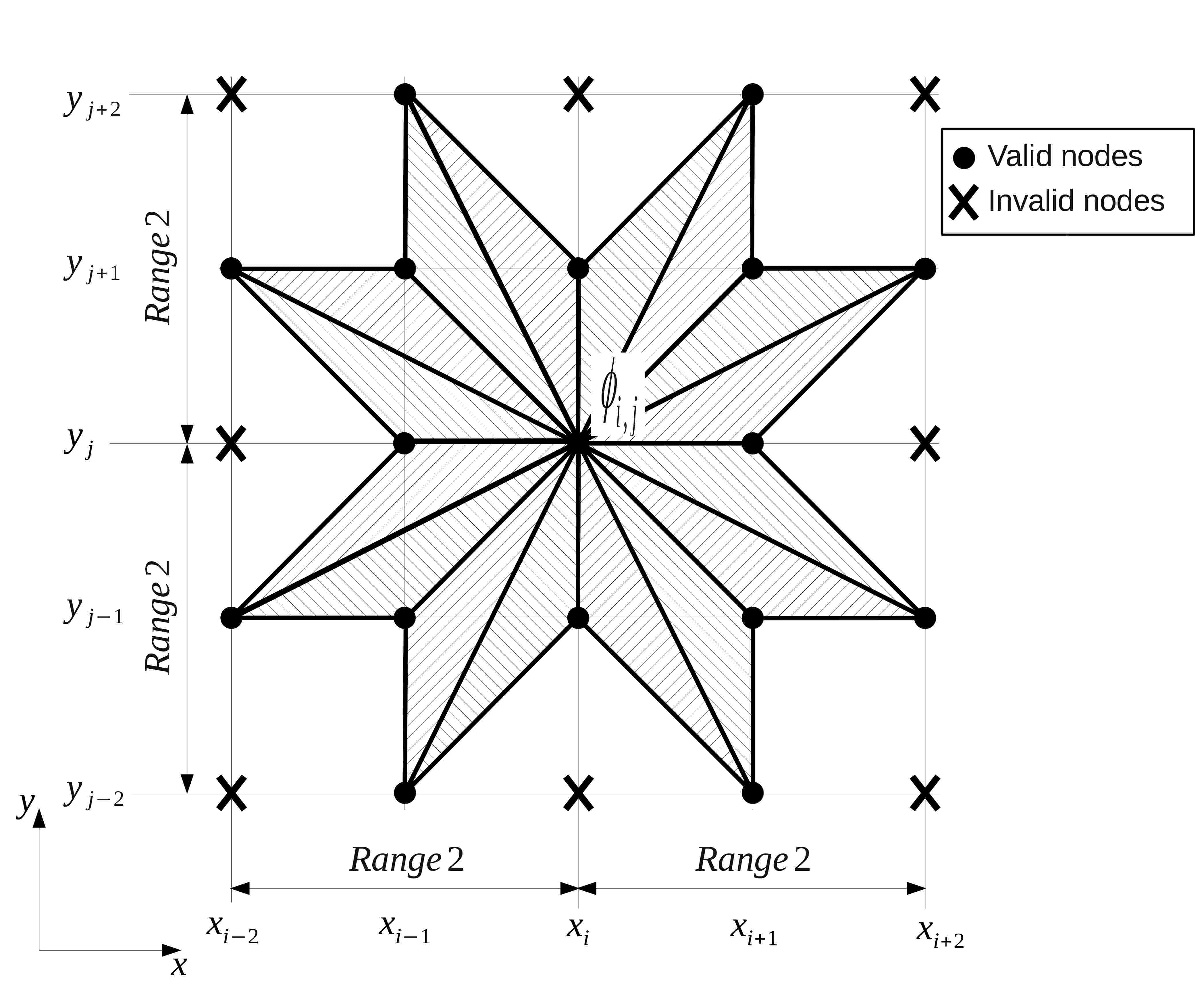}
			\caption{Shape functions for Range=2.}
		\end{subfigure}
	\end{minipage}
	\begin{subfigure}[b]{0.8\textwidth}
		\includegraphics[width=1\textwidth]{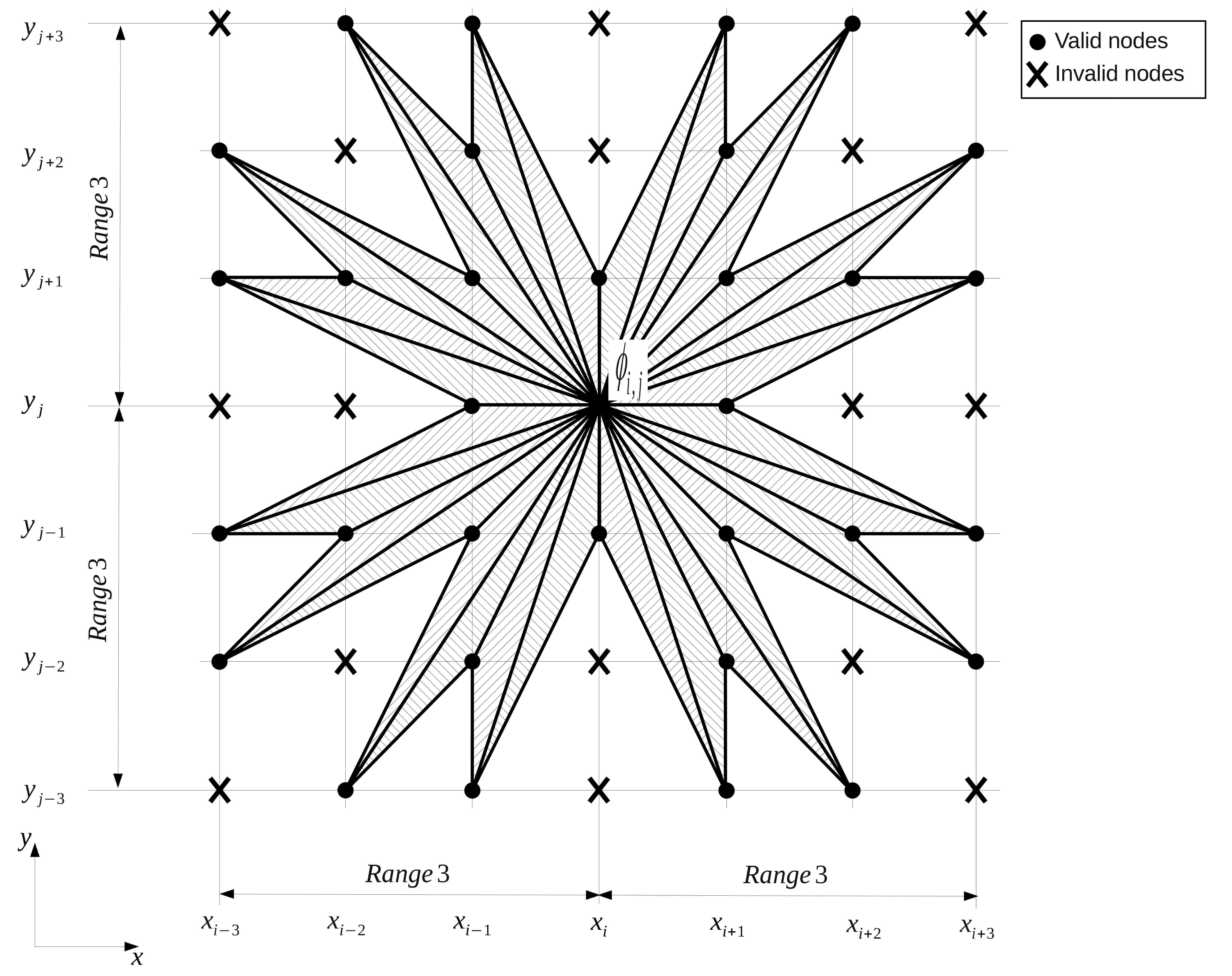}
	\caption{Shape functions for Range=3.}
	\end{subfigure}
    \caption{Shape functions according to node $\phi_{i,j}$ for ranges equal to 1 a), 2 b), and 3 c).}
    \label{Control_volumes_for_Ranges123}
\end{figure}

Two-dimensional problem is considered here. DStreaM can be extended to a three-dimensional problem straightforwardly.\\\indent
Advection equation (pure convection) can be derived from convection-diffusion equation (\ref{eq1}) when $\Gamma=0$. For simplification, density is constant ($\rho=const$):\\\indent
\begin{equation}
	\rho\frac{\partial(u \phi)}{\partial x} + \rho\frac{\partial(v \phi)}{\partial y} = 0,
	\label{eq3}
\end{equation}
DStreaM is applied on Cartesian uniform mesh for simplicity, but it can be applied as a meshless method or on an unstructured mesh. Fig. \ref{Propagation_information_cones_that_reach_node_ij_UI} shows that information from $\phi_{i-1,j}$, $\phi_{i-2,j}$ and $\phi_{i-2,j-1}$ propagates to $\phi_{i,j}$. The nodes $\phi_{i-2,j}$, $\phi_{i-1,j}$ and $\phi_{i,j}$ lie on a straight line and the information that propagates to $\phi_{i,j}$ is information from $\phi_{i-1,j}$. The information from $\phi_{i-2,j}$ does not reach $\phi_{i,j}$, because it is behind $\phi_{i-1,j}$. Therefore, the shape function interpolates nodes $\phi_{i,j}$, $\phi_{i-1,j}$ and $\phi_{i-2,j-1}$, see Fig. \ref{Propagation_information_cones_that_reach_node_ij_UI}. Convective terms shape function is a triangle with vertexes $\phi_{i,j}$, $\phi_{i-1,j}$ and $\phi_{i-2,j-1}$. It can be considered as a discrete representation of streamline through $\phi_{i,j}$, and as it is the basic element that represents the idea of the proposed approach it is included in the approach`s name. Shape function based on three nodes is:
\begin{equation}
	f^{shape}(x,y) = c_0 + c_1 x + c_2 y
	\label{eq2.0}
\end{equation}
Seeking more general form shape function nodes are denoted as $\phi_0$, $\phi_1$, and $\phi_2$ (see Fig \ref{Shape_function_and_u_average} (a)) and correspond to $\phi_{i,j}$, $\phi_{i-1,j}$, and $\phi_{i-2,j-1}$ in Fig. \ref{Propagation_information_cones_that_reach_node_ij_UI}, respectively. Coefficients $c_0$, $c_1$, and $c_2$ are determined using condition that shape function interpolates nodes $\phi_0$, $\phi_1$, and $\phi_2$ with coordinates $(x_0,y_0)$, $(x_1,y_1)$, and $(x_2,y_2)$, respectively. The coefficients are equal to:
\begin{equation}
\begin{split}
	c_0 =& \dfrac{\phi_0 (x_2 y_1 - x_1 y_2) + \phi_1 (x_0 y_2 - x_2 y_0) + \phi_2 ( x_1 y_0 - x_0 y_1)}{x_0 (y_2 - y_1) + x_1 (y_0 - y_2) + x_2 (y_1 - y_0)}\\
	c_1 =& \dfrac{\phi_0 (y_2 - y_1) + \phi_1 (y_0 - y_2) + \phi_2 (y_1 - y_0)}{x_0 (y_2 - y_1) + x_1 (y_0 - y_2) + x_2 (y_1 - y_0)}\\
	c_2 =& -\dfrac{\phi_0 (x_2 - x_1) + \phi_1 (x_0 - x_2) + \phi_2 (x_1 - x_0)}{x_0 (y_2 - y_1) + x_1 (y_0 - y_2) + x_2 (y_1 - y_0)}
	\label{eq2.1}
\end{split}
\end{equation}
\FloatBarrier
The velocities are approximated as an average value determined as follows:
\begin{equation}
\begin{split}
	u^{average} = (2 u_{i,j} + u_{i-1,j} + u_{i-2,j-1}) / 4\\
	v^{average} = (2 v_{i,j} + v_{i-1,j} + v_{i-2,j-1}) / 4
	\label{eq2}
\end{split}
\end{equation}
Fig. \ref{Shape_function_and_u_average} (b) shows $u^{average}$ geometrical representation. Control volume and shape function coincide. $u^{average}$ and $v^{average}$ are considered as a constant when integrating over control volume.  The numerical equation is derived substituting shape function (\ref{eq2.0}) and average velocities (\ref{eq2}) in (\ref{eq3}) and integrating over control volume:
\begin{equation}
\begin{split}
	&\int_{CV} \rho \frac{\partial(u \phi)}{\partial x} + \rho \frac{\partial(v \phi)}{\partial y} \,dV \approx \\
	&\rho u^{average} \int_{CV} \frac{\partial(f^{shape}(x,y))}{\partial x} \,dV + \rho v^{average} \int_{CV} \frac{\partial(f^{shape}(x,y))}{\partial y} \,dV =\\
	&\rho u^{average} c_1 V + \rho v^{average} c_2 V = 0
	\label{eq2.2}
\end{split}
\end{equation}
After substitution and simplification the numerical equation for $\phi_0$ can be expressed as:
\begin{equation}
\begin{split}
	\phi_0 = (a_1 \phi_1 + a_1 \phi_2) / a_0,
	\label{eq2.3}
\end{split}
\end{equation}
where:
\begin{equation}
\begin{split}
	a_0 =& a_1 + a_2\\
	a_1 =& u^{average} (y_0 - y_2) - v^{average} (x_0 - x_2)\\
	a_2 =& v^{average} (x_0 - x_1) - u^{average} (y_0 - y_1)
	\label{eq2.4}
\end{split}
\end{equation}
\indent
The scheme approximates values between nodes in a set range according to the approximated node ($\phi_{i,j} \equiv \phi_0$). The range applied on uniform a Cartesian mesh is a positive integer number. The nodes taken into account are on maximal distance from central node $\pm r \Delta$, where $r$ is a range, $\Delta$ is a spatial step equal to $\Delta=x_i-x_{i-1}=x_{i+1}-x_{i}=y_j-y_{j-1}=y_{j+1}-y_{j}$. Fig. \ref{Control_volumes_for_Ranges123} shows shape functions for ranges 1, 2, and 3 according to node $\phi_{i,j}$. One can see that shape functions correspond to possible discretized streamlines through node $\phi_{i,j}$ for the considered ranges. Neighbor nodes taken into account are with indexes from i-range to i+range along x-axis and from j-range to j+range along y-axis. As example when range is equal to 1 (Fig. \ref{Control_volumes_for_Ranges123} (a)), approximation scheme consider all nodes within a range $\pm \Delta$ that indexes along the x- and the y-axis are i-1,i,i+1 and j-1,j,j+1, respectively. As a result, shape functions separate space into 8 sectors. When range increases with 1, every shape function is separated into two, and a node is added. Thus, node numbers for ranges 1, 2, and 3 are 8, 16, and 32, respectively (see Fig. \ref{Control_volumes_for_Ranges123}). All angles with vertex node $\phi_{i,j}$ for range 1 are equal. For other ranges, not all angles are equal. The sector and corresponding nodes $\phi_1$ and $\phi_2$ used in approximation are selected using criteria. The selected sector has a maximum angle between average velocity vector ($\overrightarrow{V}^{average}$) and median vector $\overrightarrow{\phi_0\phi_M}$ (see Fig. \ref{Shape_function_and_u_average} (c)). In this way, the average velocity of selected shape function is maximum closer to median line $\phi_0\phi_M$, and the flow direction is from nodes $\phi_1$ and $\phi_2$ toward $\phi_0$ that corresponds to uniform information propagation of convective terms. The shape function represents the mapping of local streamline into the mesh as discrete streamline. When range increases, discrete streamline become thinner and thinner. Thinner discrete streamline tends to a line, and it is closer to a local streamline when local velocity field is uniform. When local velocity field is not uniform, thinner discrete streamlines deviate from local streamlines. Note that average velocity could be different for different sectors when velocity field is not uniform and has to be taken into account in algorithm when the sector is selected.\\\indent
\indent
Proposed approximation approach corresponds to straight local streamline that makes it simple, but some limitations occur when the velocity field is not uniform, and the range is large. DStreaM obtains solution without under/over-shoots when velocity field is uniform, and the range is large. These test cases are advection of a step profile, advection of a double-step profile, and advection of a sinusoidal profile presented further in the paper. When the local streamlines are curved, the considered approximation shows its limitations. Such approximation is acceptable for short range near the approximated node. The last considered test case is Smith and Hutton problem where local streamlines are curved, see Fig. \ref{Smith-Hutton_problem}. DStreaM obtains excellent results without under/over-shoots for ranges up to 3, see Fig. \ref{SH_contour_plot}. When the range is equal or higher than 4, obtained results have under- and over-shoots. A standard approach to prevent solution under/over-shoots and to enhance convergence is to keep coefficients $a_1$ and $a_2$ always positive, see Patel, Cross, and Markatos \cite{Patel_et_all_1988} recommendations. Here a limiter can be used, as follows:
\begin{equation}
\begin{split}
	a_1 =& max(0,u^{average} (y_0 - y_2) - v^{average} (x_0 - x_2))\\
	a_2 =& max(0,v^{average} (x_0 - x_1) - u^{average} (y_0 - y_1))
	\label{eq2.5}
\end{split}
\end{equation}
Using limiting coefficients when $a_0=a_1+a_2$ means that $\phi_0$ (see equation (\ref{eq2.3})) will be between values of $\phi_1$ and $\phi_2$, because their coefficients $a_1/a_0$ and $a_2/a_0$ are always positive, they are between 0 and 1, and their sum is equal to 1 ($a_1/a_0 + a_2/a_0 = 1$).\\\indent
Limiting coefficients prevent approximation scheme to obtain results with under/over-shoots and increases scheme ranges that obtain physically realistic results. Fig. \ref{SH_contour_plot} shows obtained fields with non-limited coefficients (\ref{eq2.4}), Fig. \ref{SH_contour_plot} (a) - (e), and limited coefficients (\ref{eq2.5}), Fig. \ref{SH_contour_plot} (g) - (k). Both obtain the same results up-to range 3. Differences occur when the range is equal to or greater than 4. For range 4 scheme with non-limited coefficients obtains a result with small under/over-shoots, while scheme with limited coefficients obtains a result without under/over-shoots, see Fig. \ref{SH_contour_plot} (e) and (j), respectively. For range 5 result obtained with non-limited coefficients is with large under- and over-shoots, Fig. \ref{SH_contour_plot} (e). The solution obtained using limited coefficients, and range 5 is without under/over-shoots, but it has unphysical jumps between some neighbor nodes, Fig. \ref{SH_contour_plot} (k). The maximal range for this spatial step that obtains correct physical results is 4 when are used limited coefficients. After all, limited coefficients make the scheme more stable and increase the scheme's application range.\\\indent
DStreaM applies to nodes with an arbitrary position near $\phi_0$ that makes is straightforwardly applicable as a meshfree method or on unstructured meshes.\\\indent
DStreaM consists of the following steps:\\\indent
1. Calculate average velocities $u^{average}$ and $v^{average}$ (\ref{eq2}).\\\indent
2. Determine upwind nodes $\phi_1$ and $\phi_2$.\\\indent
3. Calculate $\phi_0$ (\ref{eq2.3}).
\FloatBarrier
\section{Applications}
DStreaM was tested on standard advection problems: advection of a step profile, advection of a double-step profile, advection of a sinusoidal profile, and Smith and Hutton problem. Numerical equations were calculated using Gauss-Seidel iterative method with initial guess $\phi=0$. Calculation sequence of mesh nodes influences the number of iterations that each algorithm needs to reach the solution. Here were used two approaches to obtain a numerical solution with a minimal number of iterations for each scheme type. TVD schemes start every iteration from corner $i=1,j=1$. First-order upwind scheme and DStreaM change starting corner at each iteration sequentially.\\\indent
The following schemes obtained compared numerical results: first-order upwind scheme, DStreaM with ranges from 1 to 5 (DStreaM R1 to R5), TVD scheme with limiters Min-Mod, QUICK, and SUPERBEE. TVD schemes used appropriate under-relaxation coefficients to reach a solution with minimum iterations.\\\indent
DStreaM can obtain an exact solution for uniform velocity field when it uses a sufficiently large range. For presented test problems with uniform velocity field when horizontal and vertical velocities are equal ($u=v$) DStreaM obtains an exact solution for range equal to 1. These test problems are advection of a step profile, advection of a double-step profile, and advection of a sinusoidal profile. To prevent obtaining exact solution were used horizontal and vertical components of velocity equal to $u=0.8$ and $v=1$, respectively.
\subsection{Advection of a step profile}
\begin{figure}[h]
\centering
\begin{subfigure}[b]{0.37\textwidth}
	\centering
	\includegraphics[width=\textwidth]{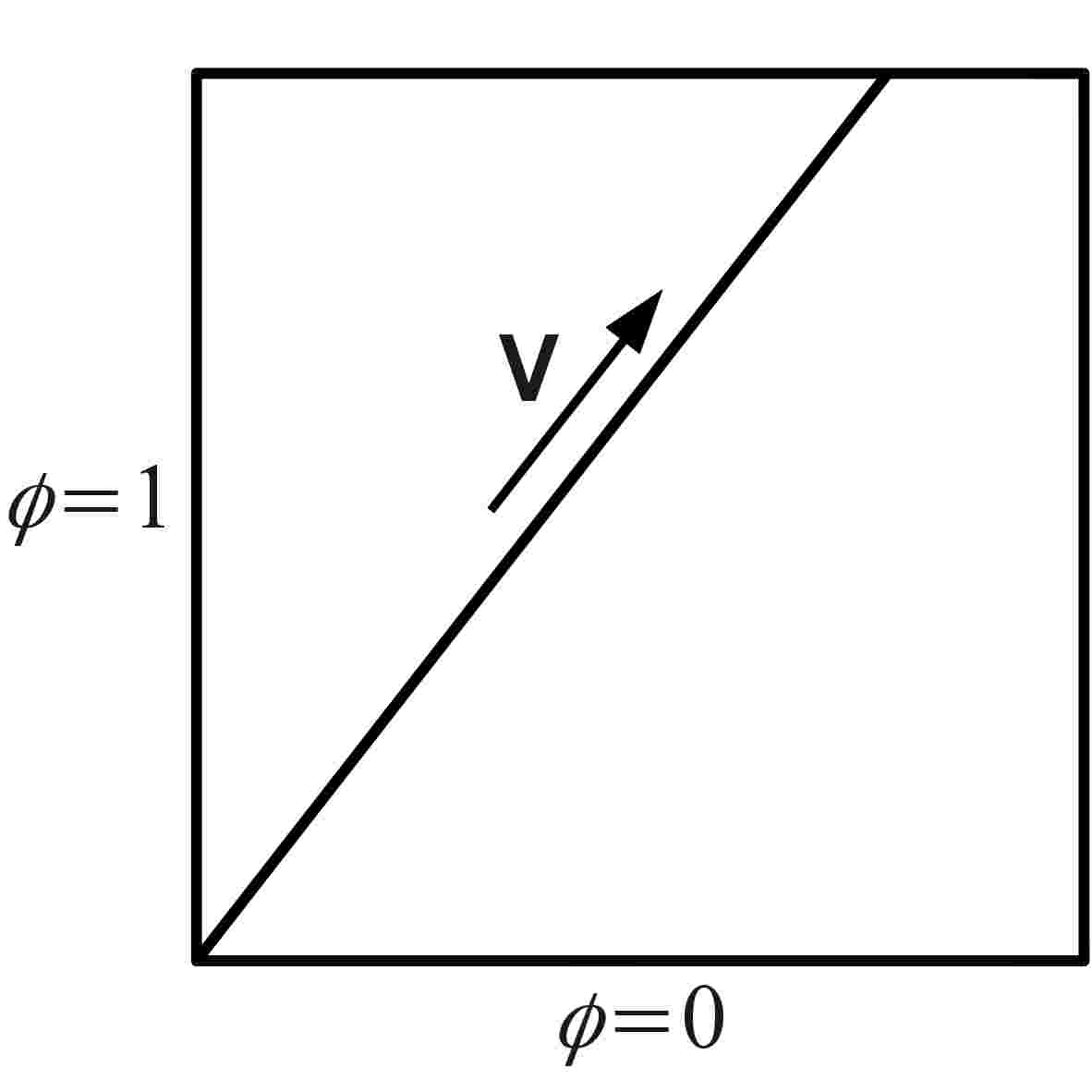}
	\caption{Physical domain\\[-8mm]}
	\label{Advection_of_a_step_profile:Physical_domain}
\end{subfigure}
\begin{subfigure}[b]{0.48\textwidth}
	\centering
	\includegraphics[width=\textwidth]{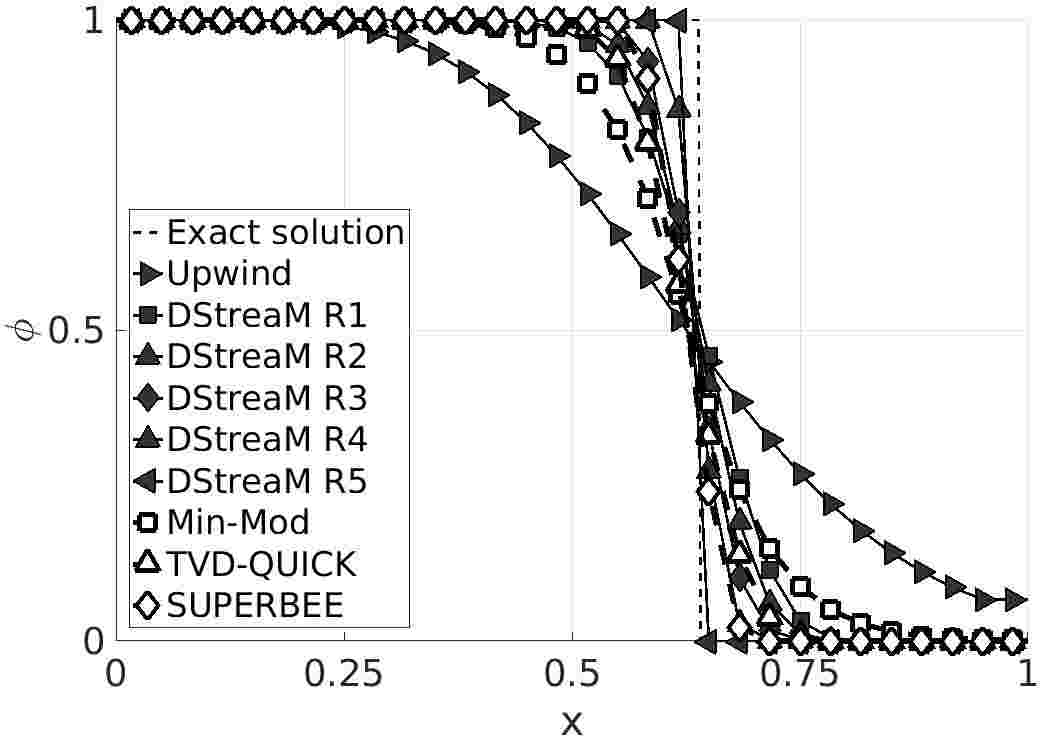}
	\caption{$\phi$ profile at $y=0.8$\\[-8mm]}
	\label{Advection_of_a_step_profile:phi_profile}
\end{subfigure}
\caption{Advection of a step profile: (a) physical domain and (b) obtained profiles at $y=0.8$.\\[-10mm]}
\label{Advection_of_a_step_profile}
\end{figure}

Horizontal and vertical components of velocity are equal to $u=0.8$ and $v=1$, respectively. The boundary conditions are rather standard for this problem, see Fig. \ref{Advection_of_a_step_profile} (a). The mesh is uniform Cartesian mesh with $30 \times 30$ nodes.\\\indent
\begin{table}[h]
\caption{Under-relaxation coefficient, number of iterations, and obtained a maximum residual of considered schemes for calculation of advection of a step profile with convergence criteria $\epsilon = 10^{-8}$.}
\centering
\begin{tabular}{lrrrr}
\toprule
Approximation scheme & $\alpha$ & Iterations & max residual \\
\midrule
Upwind     & 1.0  &     2 & 0 \\
DStreaM R1 & 1.0  &     2 & 0 \\
DStreaM R2 & 1.0  &     2 & 0 \\
DStreaM R3 & 1.0  &     2 & 0 \\
DStreaM R4 & 1.0  &     2 & 0 \\
DStreaM R5 & 1.0  &     2 & 0 \\
Min-Mod    & 0.95 &    36 & $8.1 \times 10^{-9}$ \\
QUICK(TVD) & 0.95 &    49 & $6.2 \times 10^{-9}$ \\
SUPERBEE   & 0.95 &   151 & $9.4 \times 10^{-9}$ \\
\bottomrule
\end{tabular}
\label{Advection_of_a_step_profile_Performance_30x30_10_8}
\end{table}
DStreaM with range 1 (DStreaM R1) solution is between TVD scheme with Min-Mod and QUICK limiters solutions, see Fig. \ref{Advection_of_a_step_profile} (b). DStreaM R2 solution is close to QUICK while DStreaM R3 solution is close to SUPERBEE. DStreaM R4 obtains better solution than SUPERBEE. DStreaM R5 obtains the exact solution. DStreaM does not need under-relaxation coefficient and obtains the solution for 2 iterations as upwind first-order scheme. Upwind and DStreaM need 18, 24.7, and 75.5 times fewer iterations than Min-Mod, QUICK, and SUPERBEE, see Table \ref{Advection_of_a_step_profile_Performance_30x30_10_8}.
\FloatBarrier
\subsection{Advection of a double-step profile}
\begin{figure}[h]
\centering
\begin{subfigure}[b]{0.37\textwidth}
	\centering
	\includegraphics[width=\textwidth]{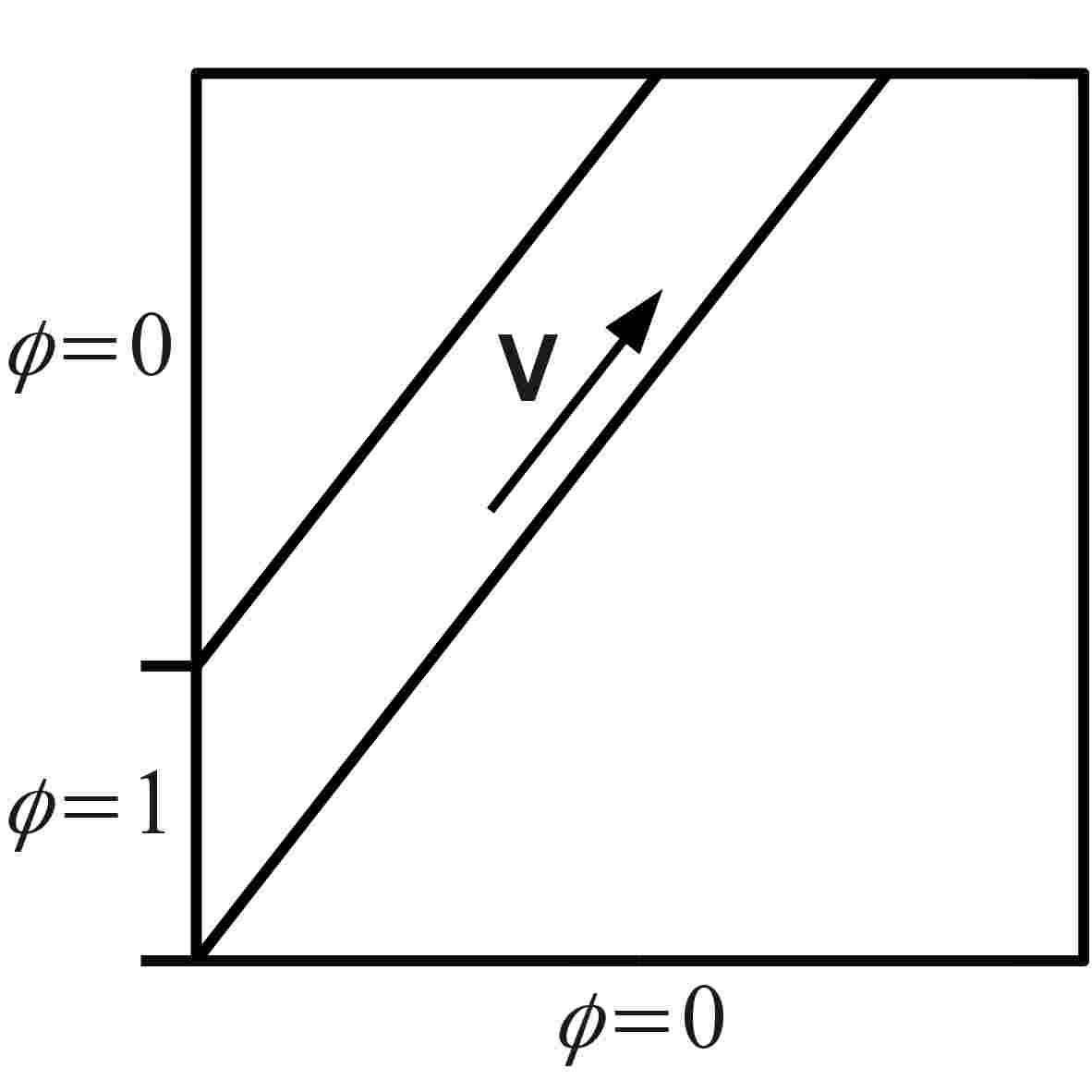}\\[-5mm]
	\caption{Physical domain\\[-7mm]}
	\label{Advection_of_a_double-step_profile:Physical_domain}
\end{subfigure}
\begin{subfigure}[b]{0.48\textwidth}
	\centering
	\includegraphics[width=\textwidth]{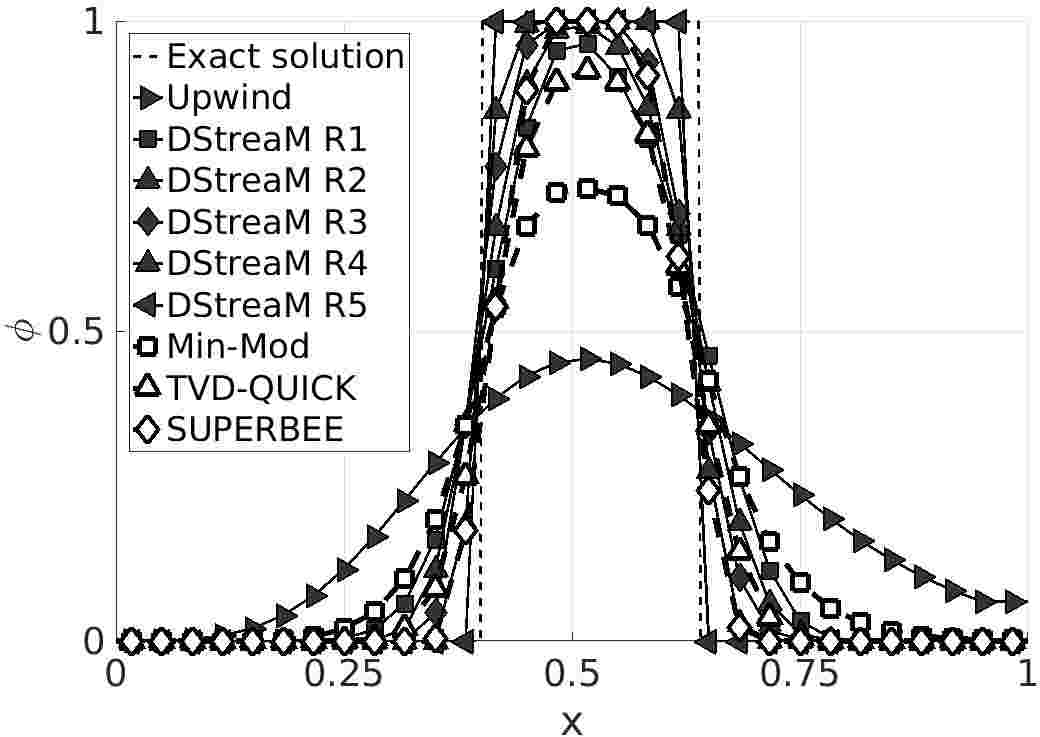}\\[-5mm]
	\caption{$\phi$ profile at $y=0.8$\\[-7mm]}
	\label{Advection_of_a_double-step_profile:phi_profile}
\end{subfigure}
\caption{Advection of a double-step profile: (a) physical domain and (b) obtained profiles at $y=0.8$.\\[-10mm]}
\label{Advection_of_a_double-step_profile}
\end{figure}

Horizontal and vertical components of velocity are equal to $u=0.8$ and $v=1$, respectively. The boundary conditions are rather standard for this problem, see Fig. \ref{Advection_of_a_double-step_profile} (a). The mesh is uniform Cartesian mesh with $30 \times 30$ nodes.\\\indent
\begin{table}[h]
\caption{Under-relaxation coefficient, number of iterations, and obtained a maximum residual of considered schemes for calculation of advection of a double-step profile with convergence criteria $\epsilon = 10^{-8}$.}
\centering
\begin{tabular}{lrrrr}
\toprule
Approximation scheme & $\alpha$ & Iterations & max residual \\
\midrule
Upwind     & 1.0  &     2 & 0 \\
DStreaM R1 & 1.0  &     2 & 0 \\
DStreaM R2 & 1.0  &     2 & 0 \\
DStreaM R3 & 1.0  &     2 & 0 \\
DStreaM R4 & 1.0  &     2 & 0 \\
DStreaM R5 & 1.0  &     2 & 0 \\
Min-Mod    & 1.0  &    30 & $6.7 \times 10^{-9}$ \\
QUICK(TVD) & 0.95 &    55 & $5.6 \times 10^{-9}$ \\
SUPERBEE   & 0.9  &   187 & $9.9 \times 10^{-9}$ \\
\bottomrule
\end{tabular}
\label{Advection_of_a_double-step_profile_Performance_30x30_10_8}
\end{table}
DStreaM R1 solution is between solutions obtained by TVD scheme with QUICK and SUPERBEE limiters, see Fig. \ref{Advection_of_a_double-step_profile} (b). DStreaM R3 solution is close to SUPERBEE. DStreaM R4 obtains better solution than SUPERBEE. DStreaM R5 obtains the exact solution. DStreaM does not need under-relaxation coefficient and obtains the solution for 2 iterations as upwind first-order scheme. Upwind and DStreaM need 15, 27.5, and 93.5 times fewer iterations than Min-Mod, QUICK, and SUPERBEE, see Table \ref{Advection_of_a_double-step_profile_Performance_30x30_10_8}.
\FloatBarrier
\subsection{Advection of a sinusoidal profile}
\begin{figure}[h]
\centering
\begin{subfigure}[b]{0.51\textwidth}
	\centering
	\includegraphics[width=\textwidth]{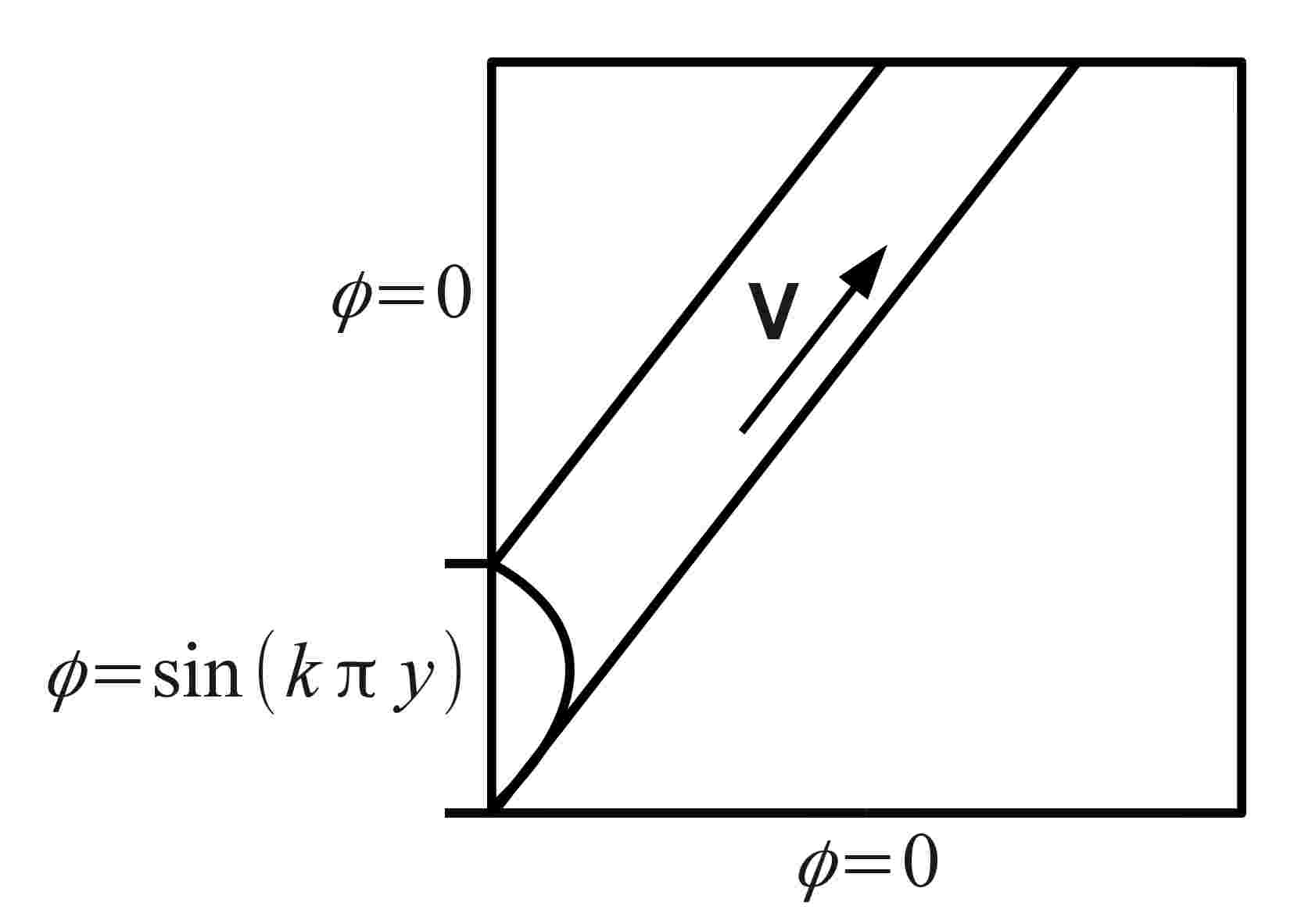}
	\caption{Physical domain\\[-8mm]}
	\label{Advection_of_a_sinusoidal_profile:Physical_domain}
\end{subfigure}
\begin{subfigure}[b]{0.48\textwidth}
	\centering
	\includegraphics[width=\textwidth]{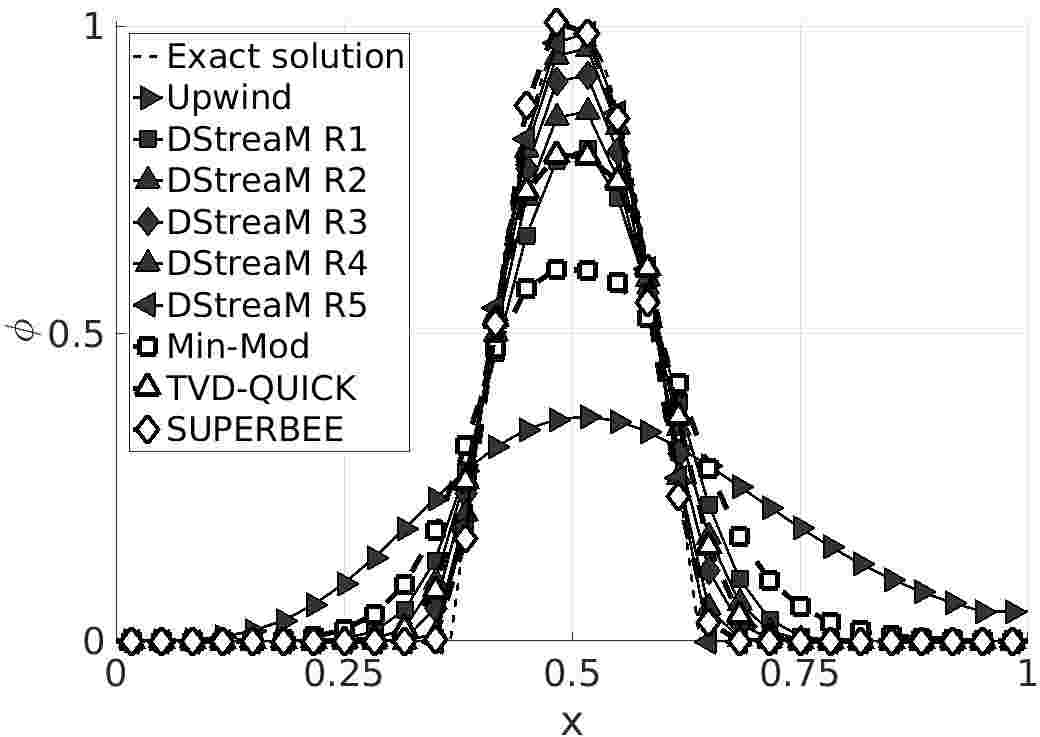}
	\caption{$\phi$ profile at $y=0.8$\\[-8mm]}
	\label{Advection_of_a_sinusoidal_profile:phi_profile}
\end{subfigure}
\caption{Advection of a sinusoidal profile: (a) physical domain and (b) obtained profiles at $y=0.8$.\\[-8mm]}
\label{Advection_of_a_sinusoidal_profile}
\end{figure}

Horizontal and vertical components of velocity are equal to $u=0.8$ and $v=1$, respectively. The boundary conditions are rather standard for this problem, see Fig. \ref{Advection_of_a_sinusoidal_profile} (a). The mesh is uniform Cartesian mesh with $30 \times 30$ nodes.\\\indent
\begin{table}[h]
\caption{Under-relaxation coefficient, number of iterations, and obtained a maximum residual of considered schemes for calculation of advection of a sinusoidal profile with convergence criteria $\epsilon = 10^{-8}$.}
\centering
\begin{tabular}{lrrrr}
\toprule
Approximation scheme & $\alpha$ & Iterations & max residual \\
\midrule
Upwind     & 1.0  &     2 & 0 \\
DStreaM R1 & 1.0  &     2 & 0 \\
DStreaM R2 & 1.0  &     2 & 0 \\
DStreaM R3 & 1.0  &     2 & 0 \\
DStreaM R4 & 1.0  &     2 & 0 \\
DStreaM R5 & 1.0  &     2 & 0 \\
Min-Mod    & 0.95 &    30 & $6.6 \times 10^{-9}$ \\
QUICK(TVD) & 0.9  &    51 & $7.9 \times 10^{-9}$ \\
SUPERBEE   & 0.85 &    65 & $8.3 \times 10^{-9}$ \\
\bottomrule
\end{tabular}
\label{Advection_of_a_sinusoidal_profile_Performance_30x30_10_8}
\end{table}
DStreaM R1 solution corresponds to QUICK's solution, see Fig. \ref{Advection_of_a_double-step_profile} (b). DStreaM R2, R3, and R4 solutions are closer to the exact solution than QUICK's and SUPERBEE's solutions. DStreaM R5 obtains the exact solution. DStreaM scheme does not need under-relaxation coefficient and obtains a solution for 2 iterations as upwind first-order scheme. Upwind and DStreaM need 15, 25.5, and 32.5 times fewer iterations than Min-Mod, QUICK, and SUPERBEE, see Table \ref{Advection_of_a_sinusoidal_profile_Performance_30x30_10_8}.
\FloatBarrier
\subsection{Smith and Hutton problem}
\begin{figure}[htb!]
	\centering \ \\[-4mm]
	\includegraphics[height=0.4\textwidth]{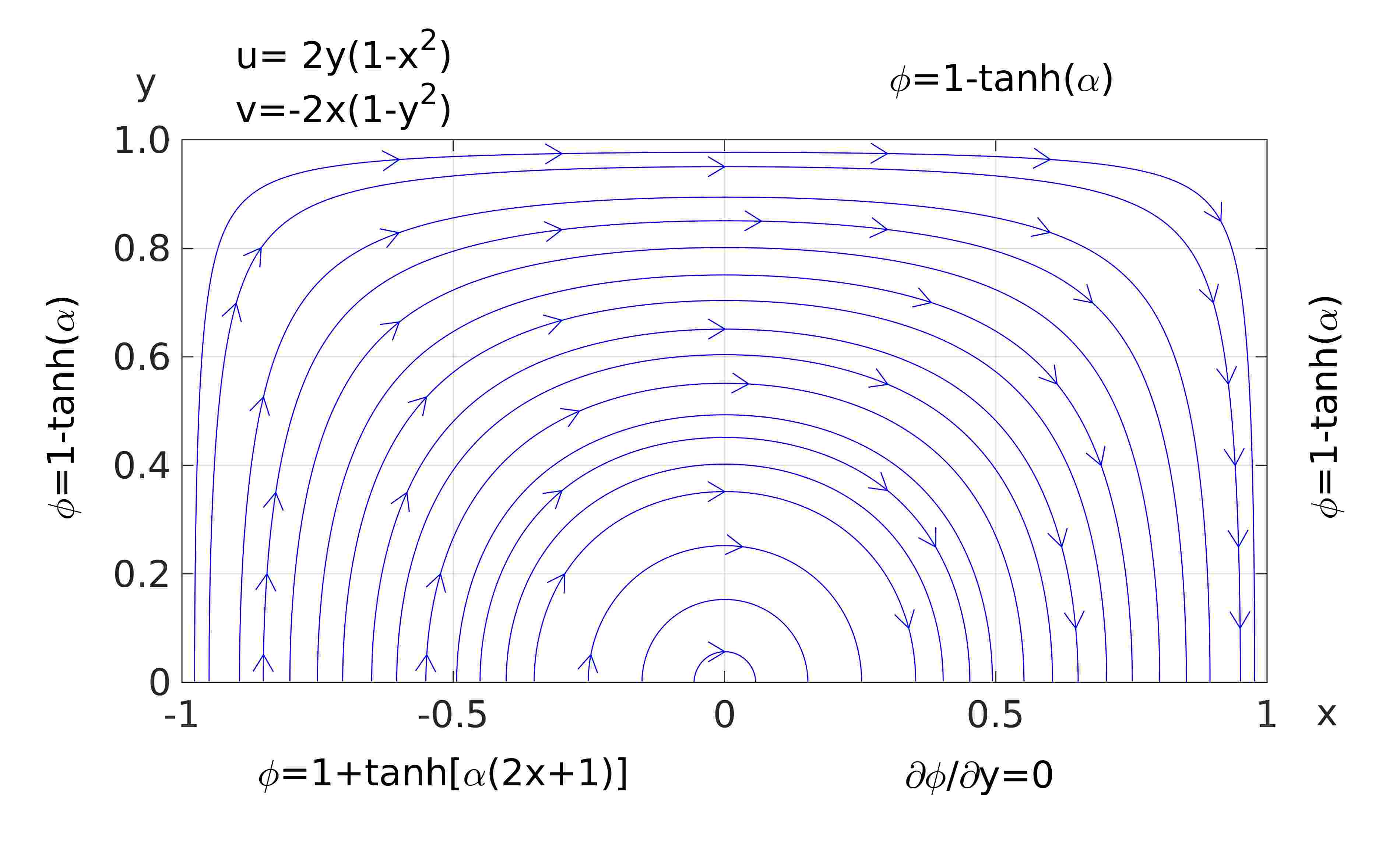}\\[-8mm]
	\caption{Smith and Hutton problem: boundary conditions, velocity equations, and streamlines.\\[-5mm]}
	\label{Smith-Hutton_problem}
\end{figure}

\begin{figure}[htbp!]
	\centering
	\tiny
	\begin{minipage}{1.2\textwidth}
	\raggedright
		\begin{minipage}{0.3\textwidth}
			\ \\[30mm]
			\begin{subfigure}[b]{1.0\textwidth}
				\rotatebox{90}{\parbox{2.5cm}{\raggedright (a) DStreaM R1,\\non-limited coeff. (\ref{eq2.4})}}
				\includegraphics[trim = 134 0 350 280,clip=true, keepaspectratio=true,width=0.9\textwidth]{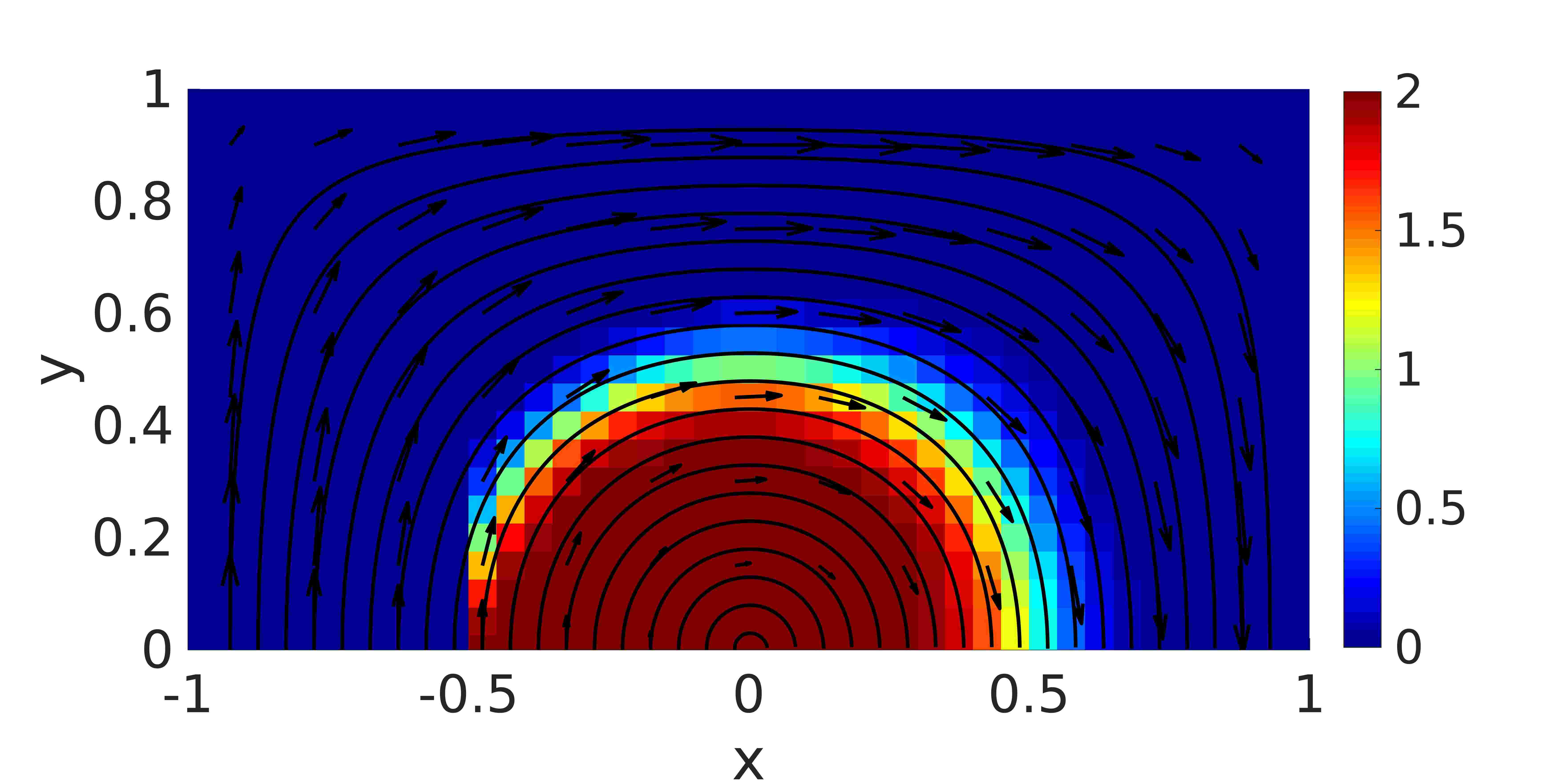}
			\end{subfigure}
			\begin{subfigure}[b]{1.0\textwidth}
				\rotatebox{90}{\parbox{2.5cm}{\raggedright (b) DStreaM R2,\\non-limited coeff. (\ref{eq2.4})}}
				\includegraphics[trim = 134 0 350 280,clip=true, keepaspectratio=true,width=0.9\textwidth]{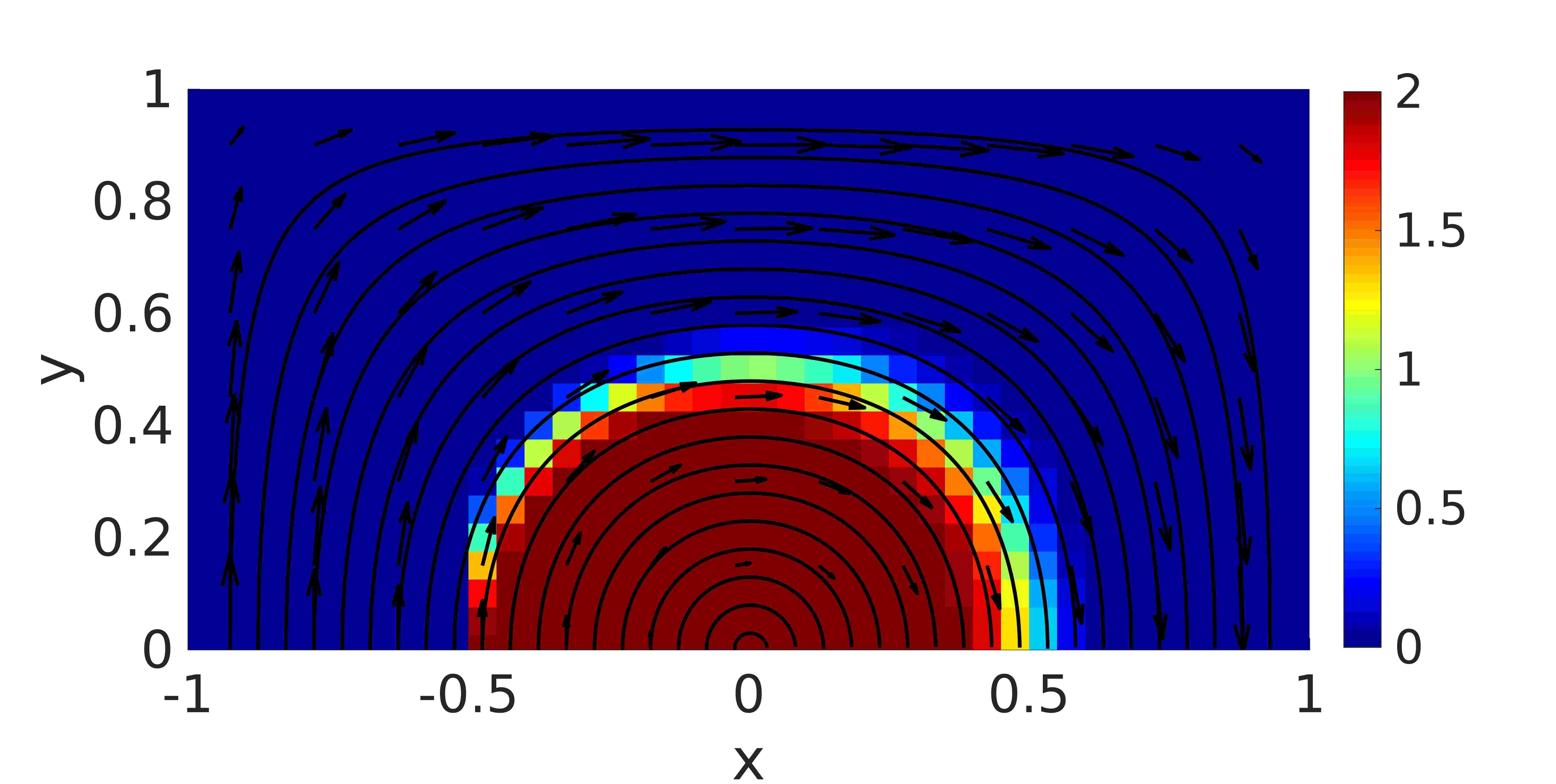}
			\end{subfigure}
			\begin{subfigure}[b]{1.0\textwidth}
				\rotatebox{90}{\parbox{2.5cm}{\raggedright (c) DStreaM R3,\\non-limited coeff. (\ref{eq2.4})}}
				\includegraphics[trim = 134 0 350 280,clip=true, keepaspectratio=true,width=0.9\textwidth]{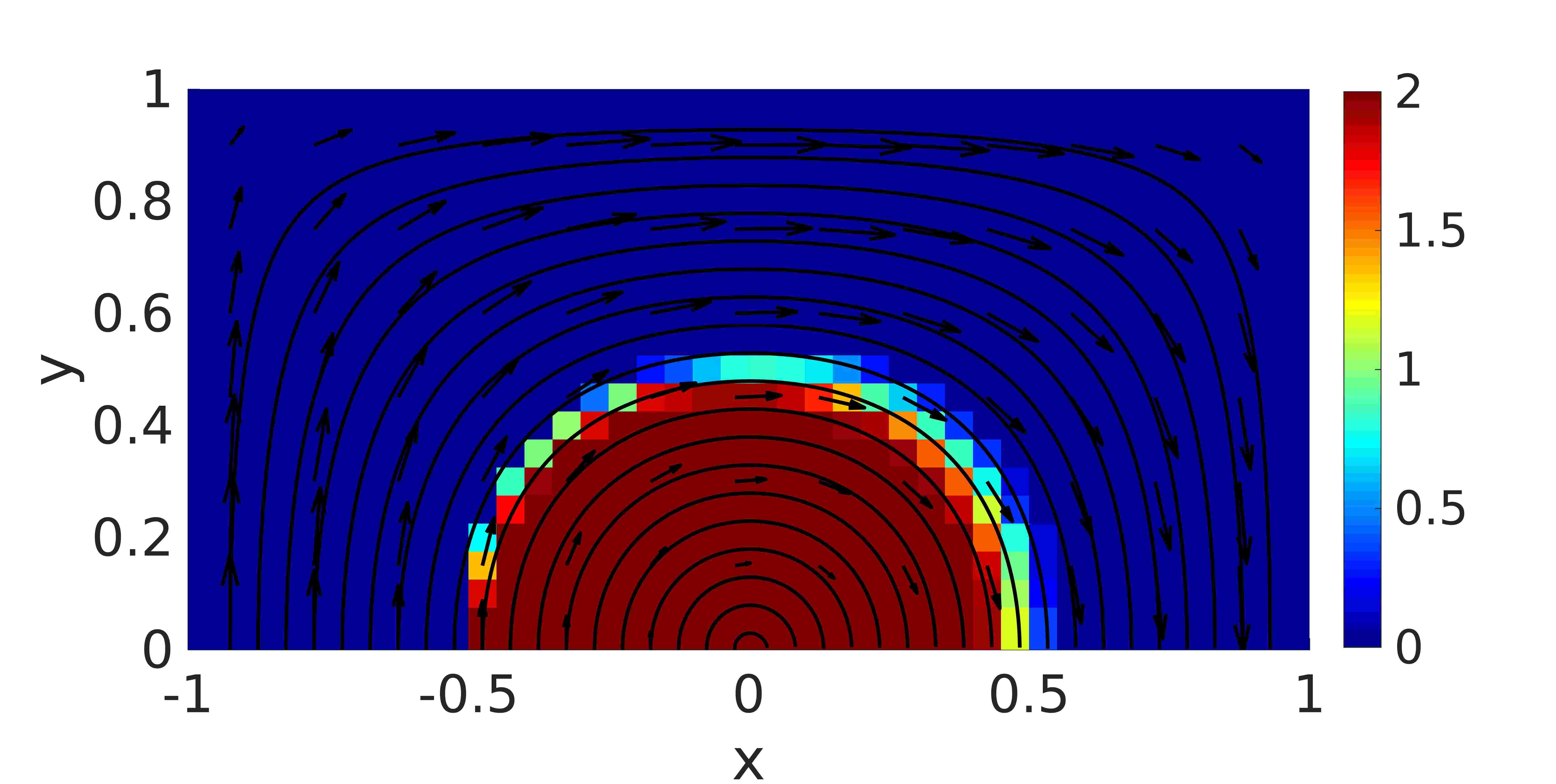}
			\end{subfigure}
			\begin{subfigure}[b]{1.0\textwidth}
				\rotatebox{90}{\parbox{2.5cm}{\raggedright (d) DStreaM R4,\\non-limited coeff. (\ref{eq2.4})}}
				\includegraphics[trim = 134 0 350 280,clip=true, keepaspectratio=true,width=0.9\textwidth]{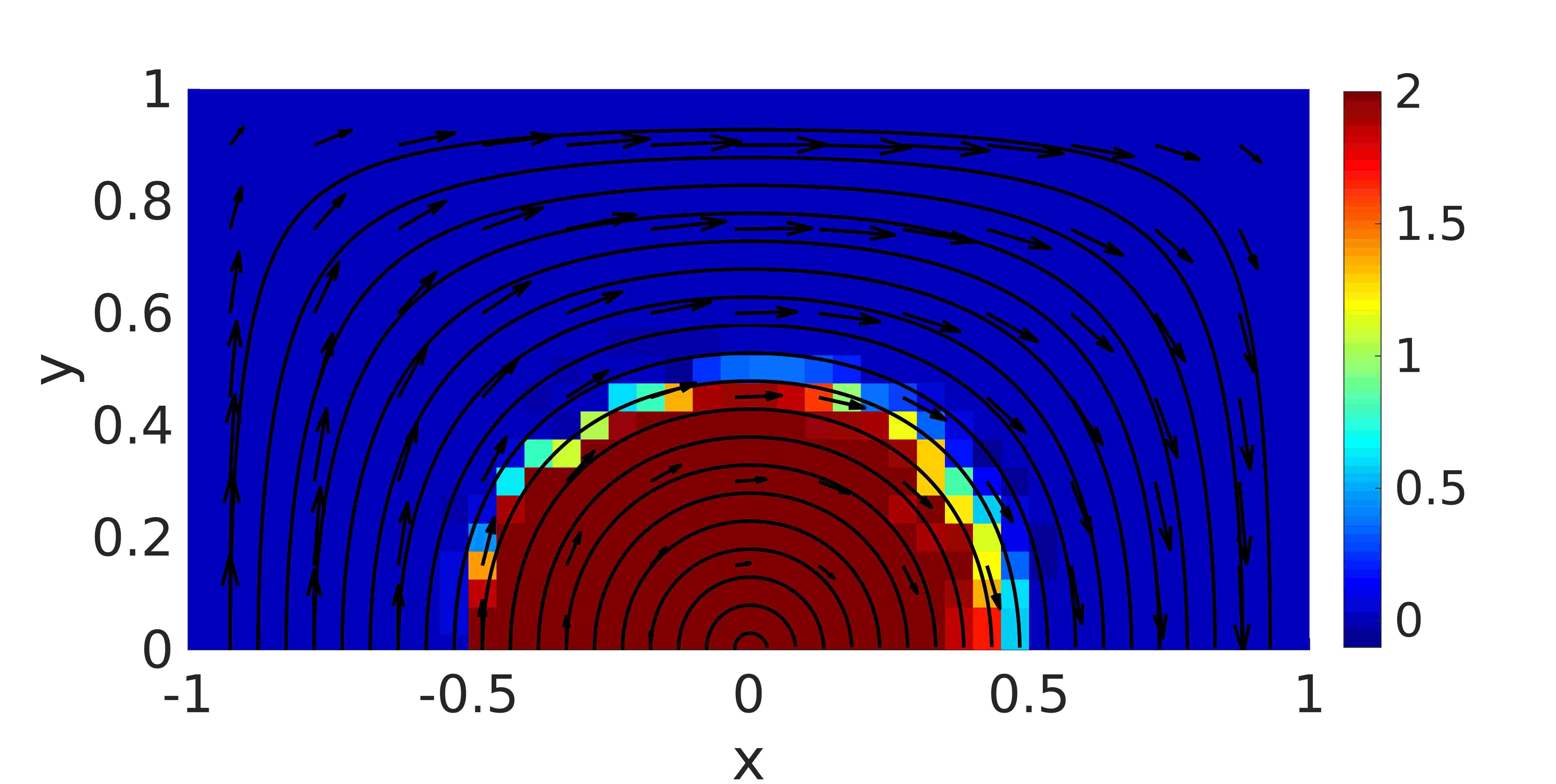}
			\end{subfigure}
			\begin{subfigure}[b]{1.0\textwidth}
				\rotatebox{90}{\parbox{2.5cm}{\raggedright (e) DStreaM R5,\\non-limited coeff. (\ref{eq2.4})}}
				\includegraphics[trim = 134 0 350 280,clip=true, keepaspectratio=true,width=0.9\textwidth]{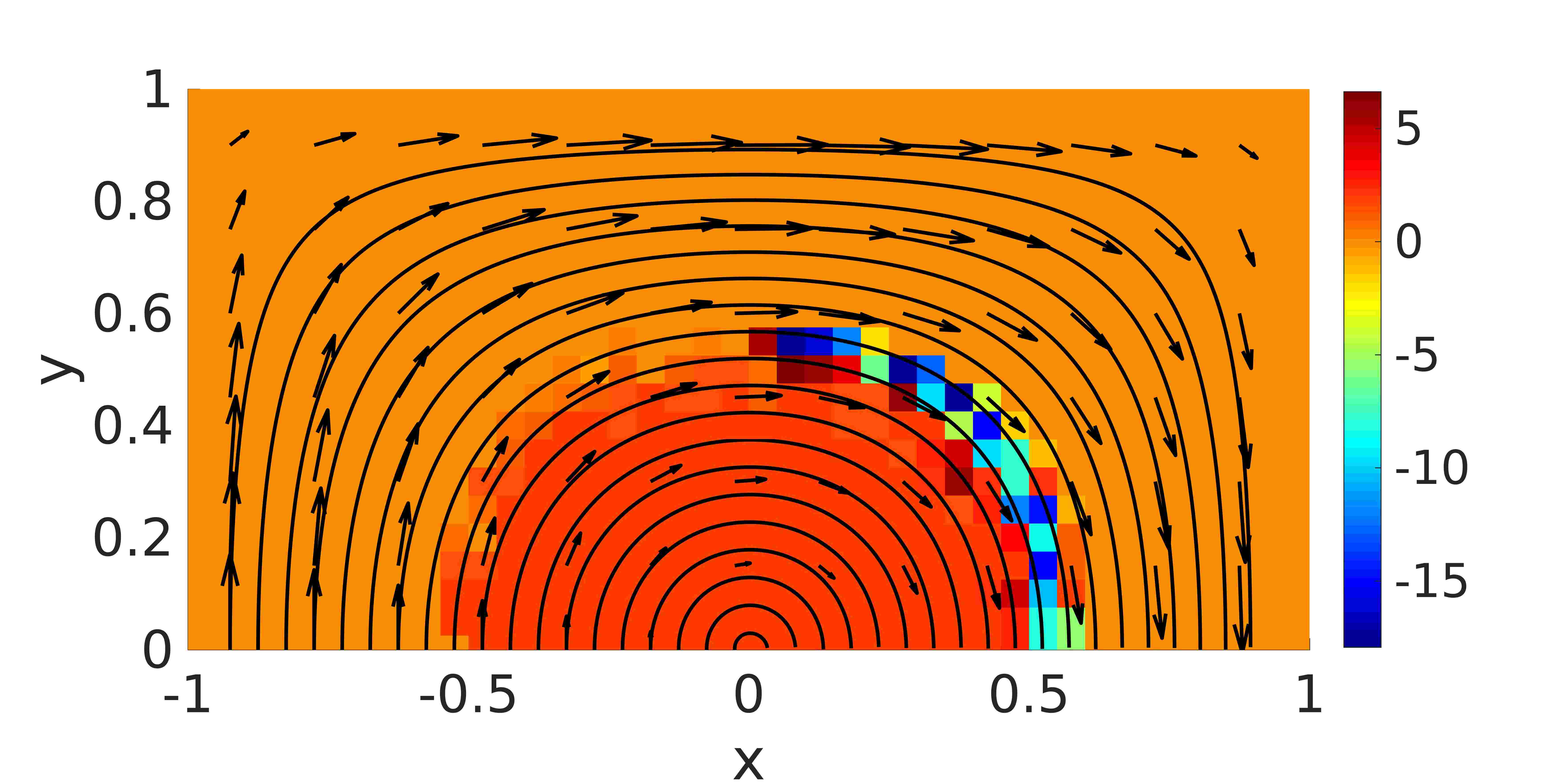}
			\end{subfigure}
		\end{minipage}
		\begin{minipage}{0.3\textwidth}
			\begin{subfigure}[b]{1.0\textwidth}
				\rotatebox{90}{\parbox{2.5cm}{\raggedright (f) Upwind\\ \ }}
				\includegraphics[trim = 134 0 350 280,clip=true, keepaspectratio=true,width=0.9\textwidth]{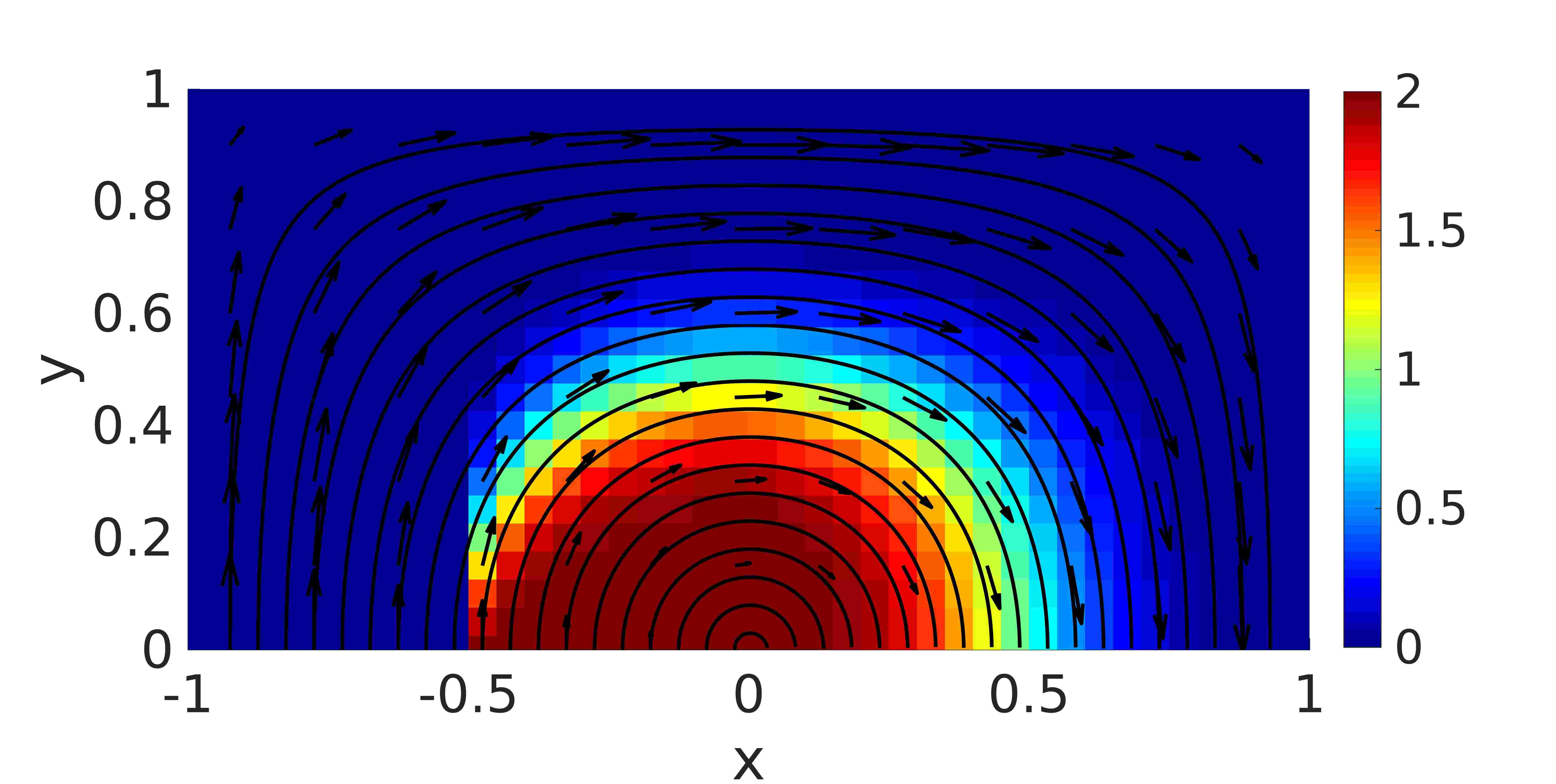}
			\end{subfigure}
			\begin{subfigure}[b]{1.0\textwidth}
				\rotatebox{90}{\parbox{2.5cm}{\raggedright (g) DStreaM R1,\\limited coeff. (\ref{eq2.5})}}
				\includegraphics[trim = 134 0 350 280,clip=true, keepaspectratio=true,width=0.9\textwidth]{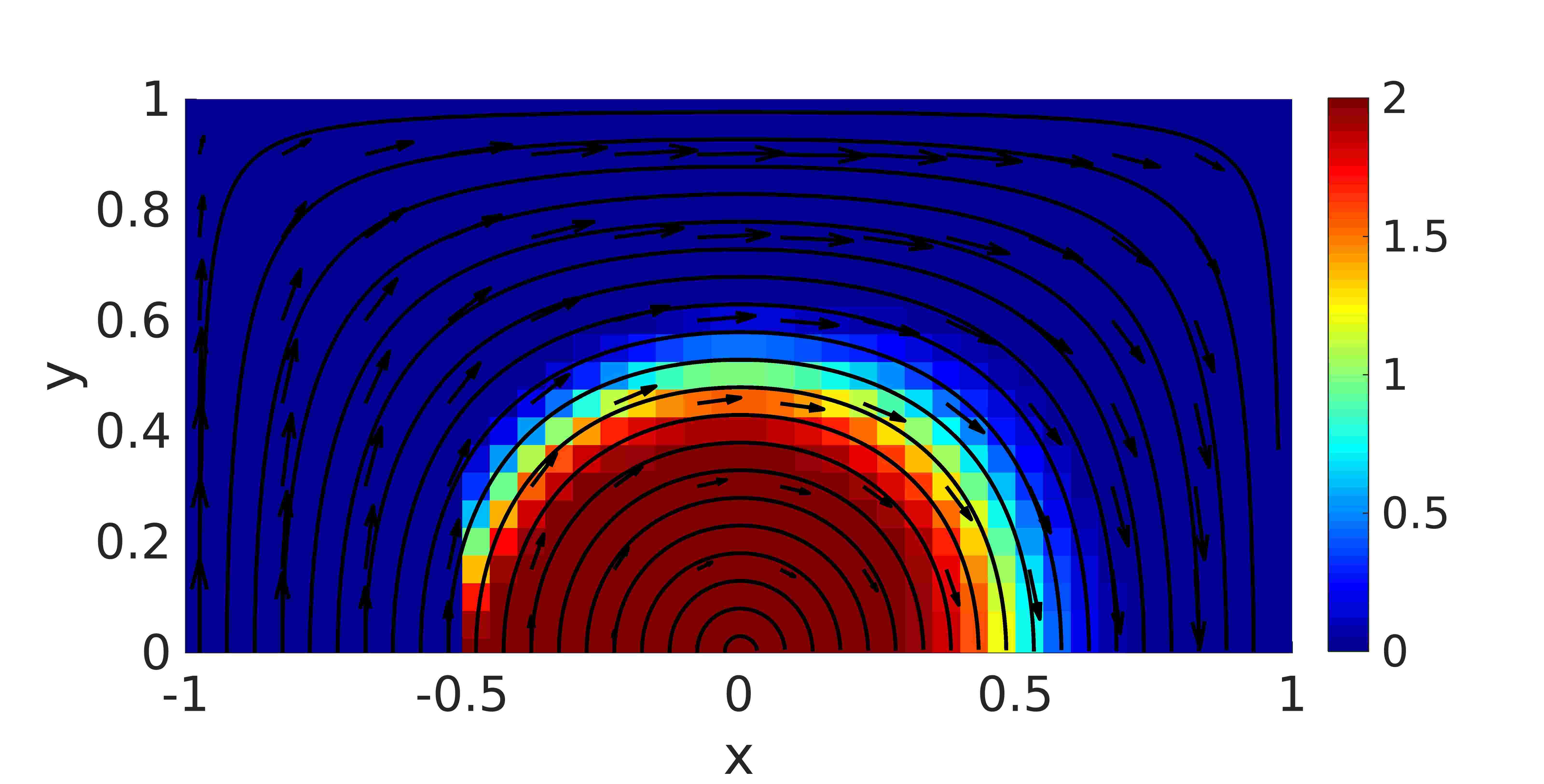}
			\end{subfigure}
			\begin{subfigure}[b]{1.0\textwidth}
				\rotatebox{90}{\parbox{2.5cm}{\raggedright (h) DStreaM R2,\\limited coeff. (\ref{eq2.5})}}
				\includegraphics[trim = 134 0 350 280,clip=true, keepaspectratio=true,width=0.9\textwidth]{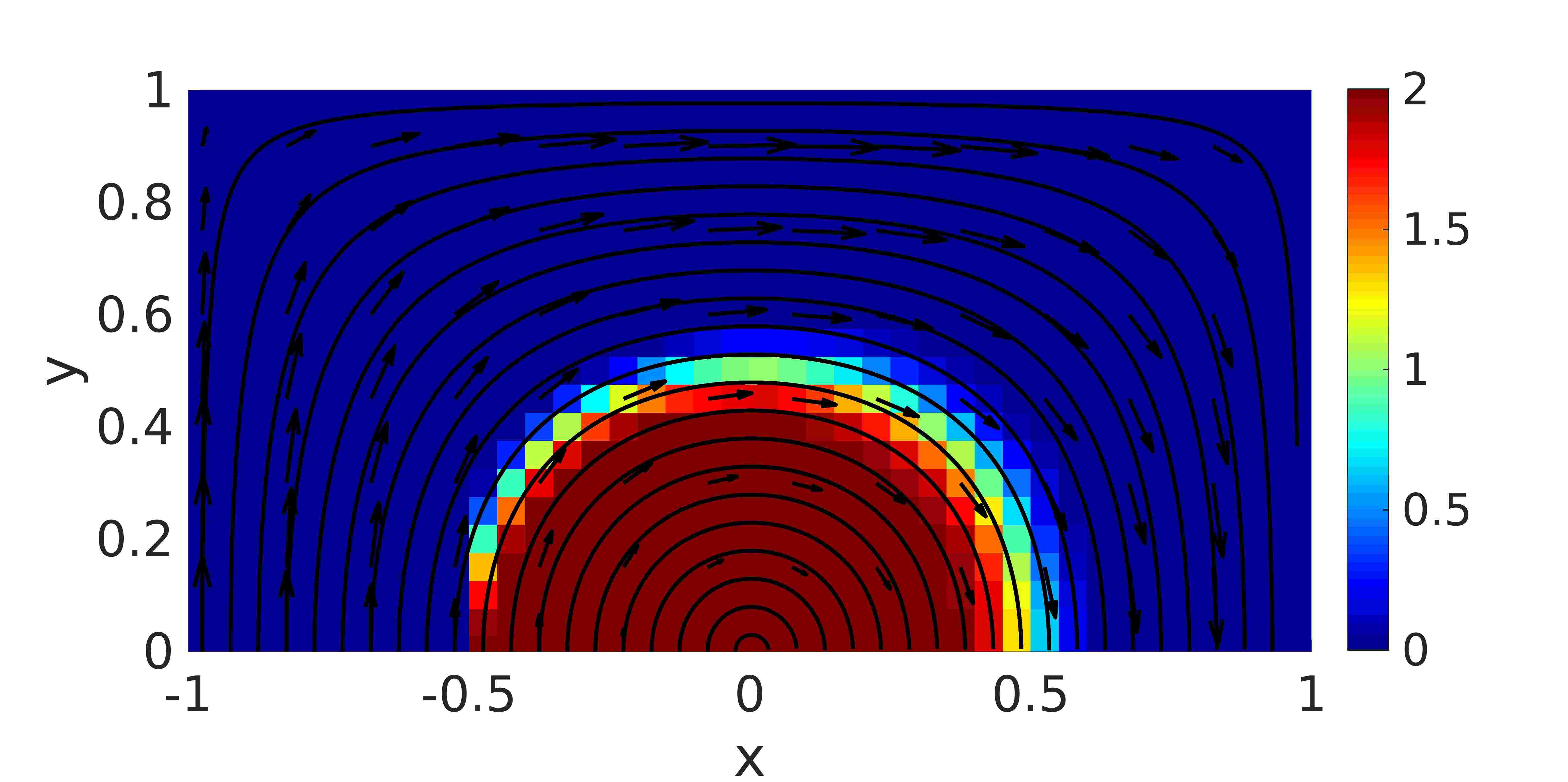}
			\end{subfigure}
			\begin{subfigure}[b]{1.0\textwidth}
				\rotatebox{90}{\parbox{2.5cm}{\raggedright (i) DStreaM R3,\\limited coeff. (\ref{eq2.5})}}
				\includegraphics[trim = 134 0 350 280,clip=true, keepaspectratio=true,width=0.9\textwidth]{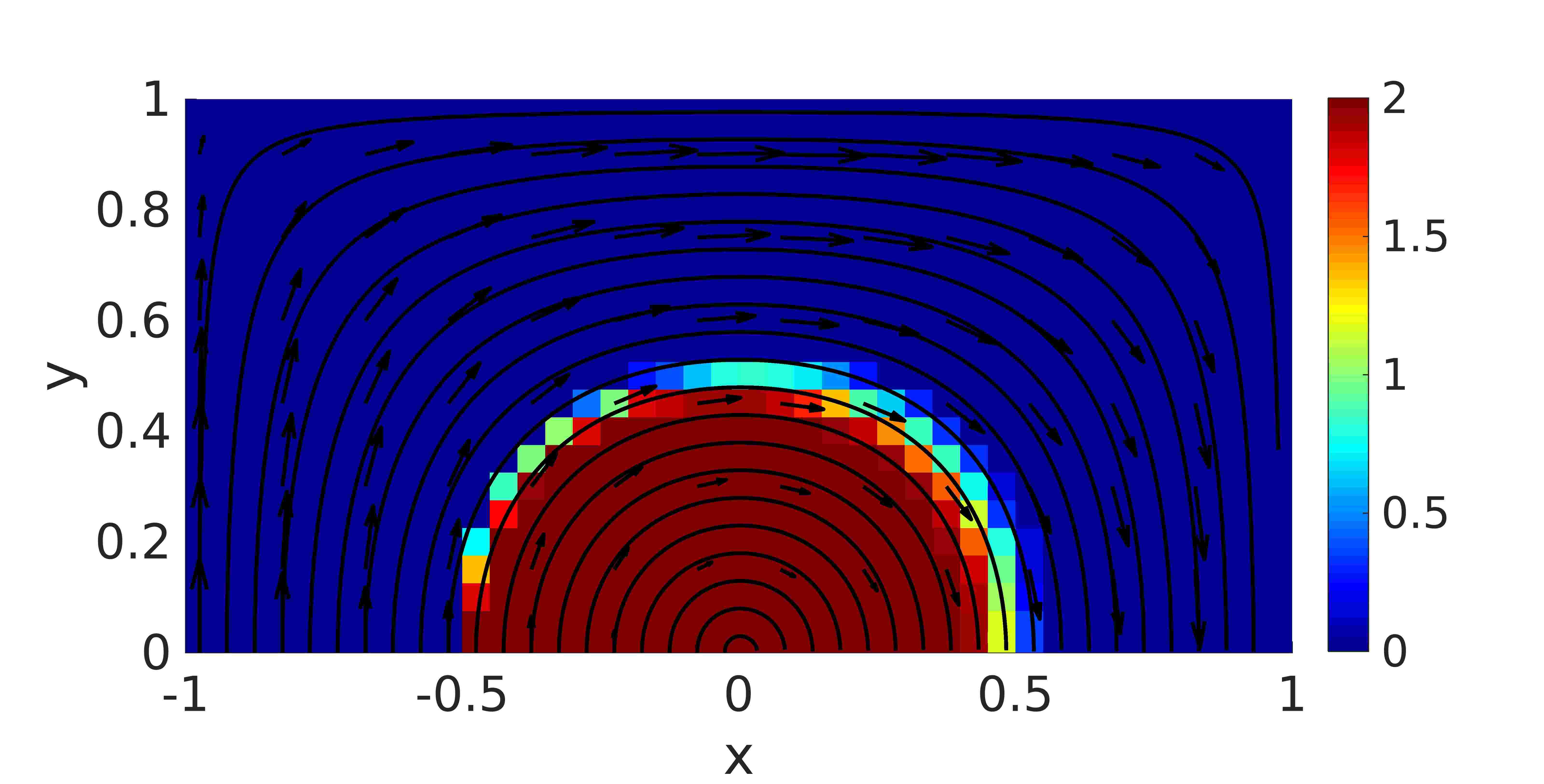}
			\end{subfigure}
			\begin{subfigure}[b]{1.0\textwidth}
				\rotatebox{90}{\parbox{2.5cm}{\raggedright (j) DStreaM R4,\\limited coeff. (\ref{eq2.5})}}
				\includegraphics[trim = 134 0 350 280,clip=true, keepaspectratio=true,width=0.9\textwidth]{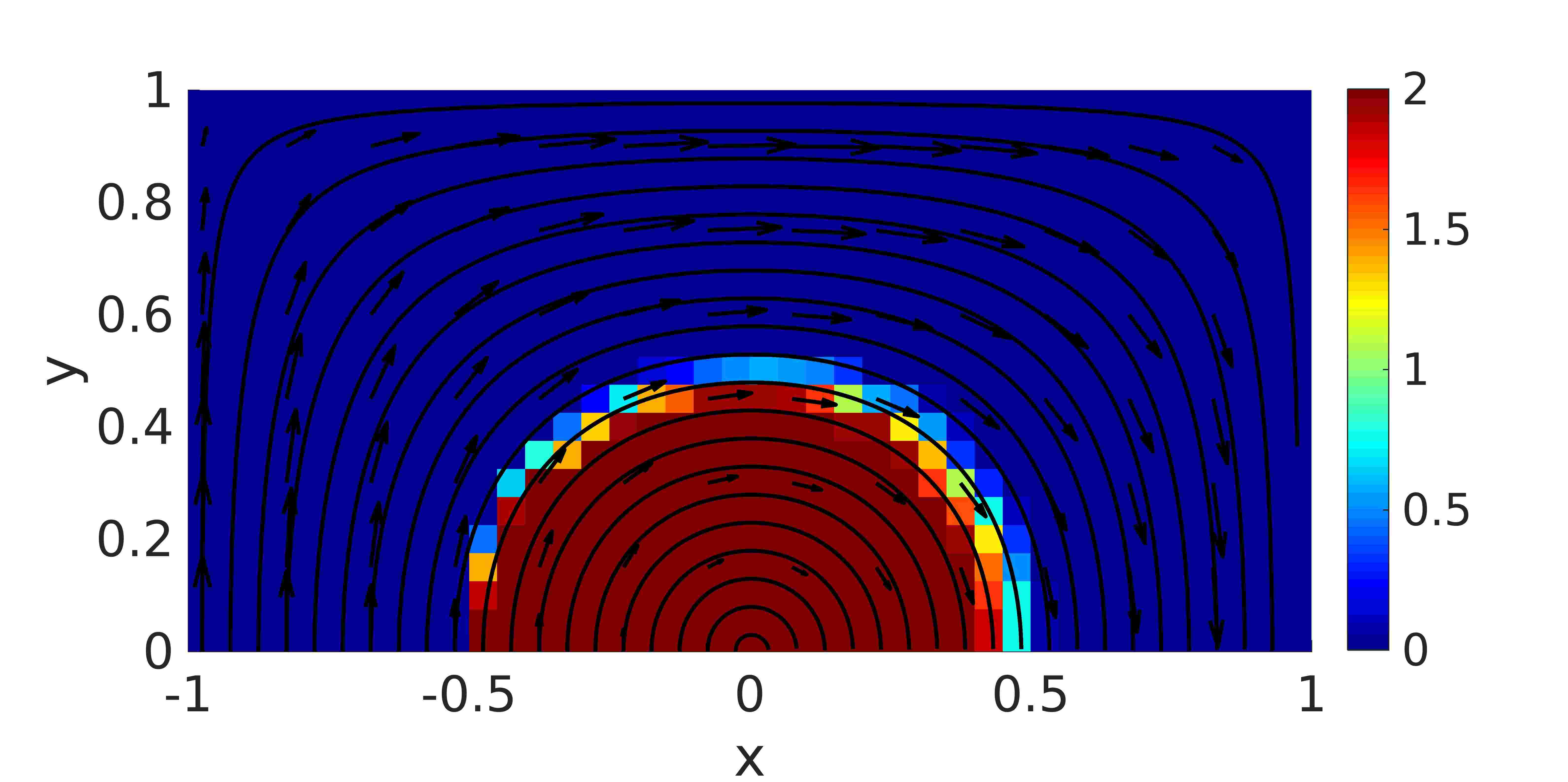}
			\end{subfigure}
			\begin{subfigure}[b]{1.0\textwidth}
				\rotatebox{90}{\parbox{2.5cm}{\raggedright (k) DStreaM R5,\\limited coeff. (\ref{eq2.5})}}
				\includegraphics[trim = 134 0 350 280,clip=true, keepaspectratio=true,width=0.9\textwidth]{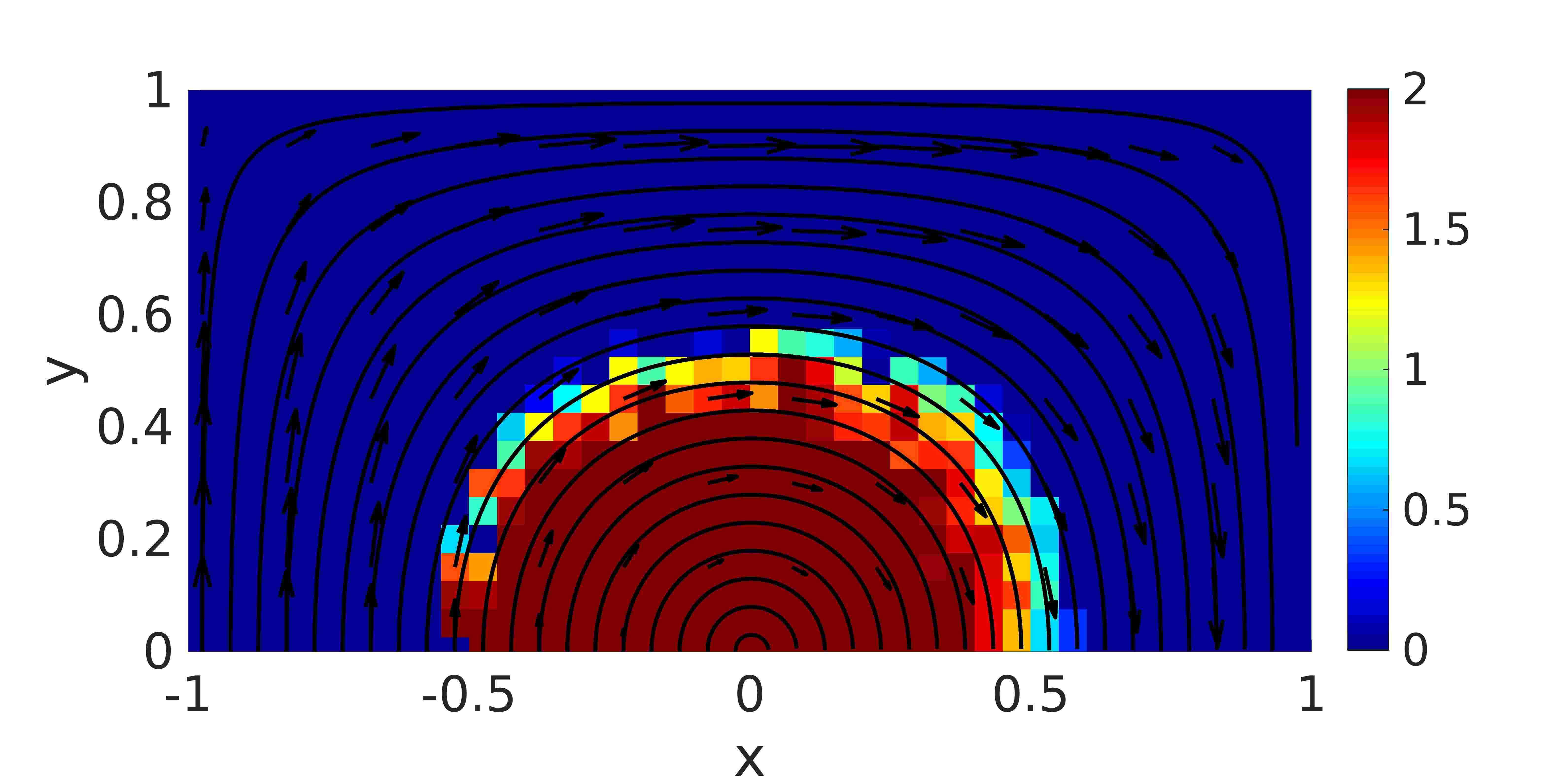}
			\end{subfigure}
		\end{minipage}
		\begin{minipage}{0.3\textwidth}
			\ \\[10mm]
			\begin{subfigure}[b]{1.0\textwidth}
				\rotatebox{90}{(l) Min-Mod}
				\includegraphics[trim = 134 0 350 280,clip=true, keepaspectratio=true,width=0.9\textwidth]{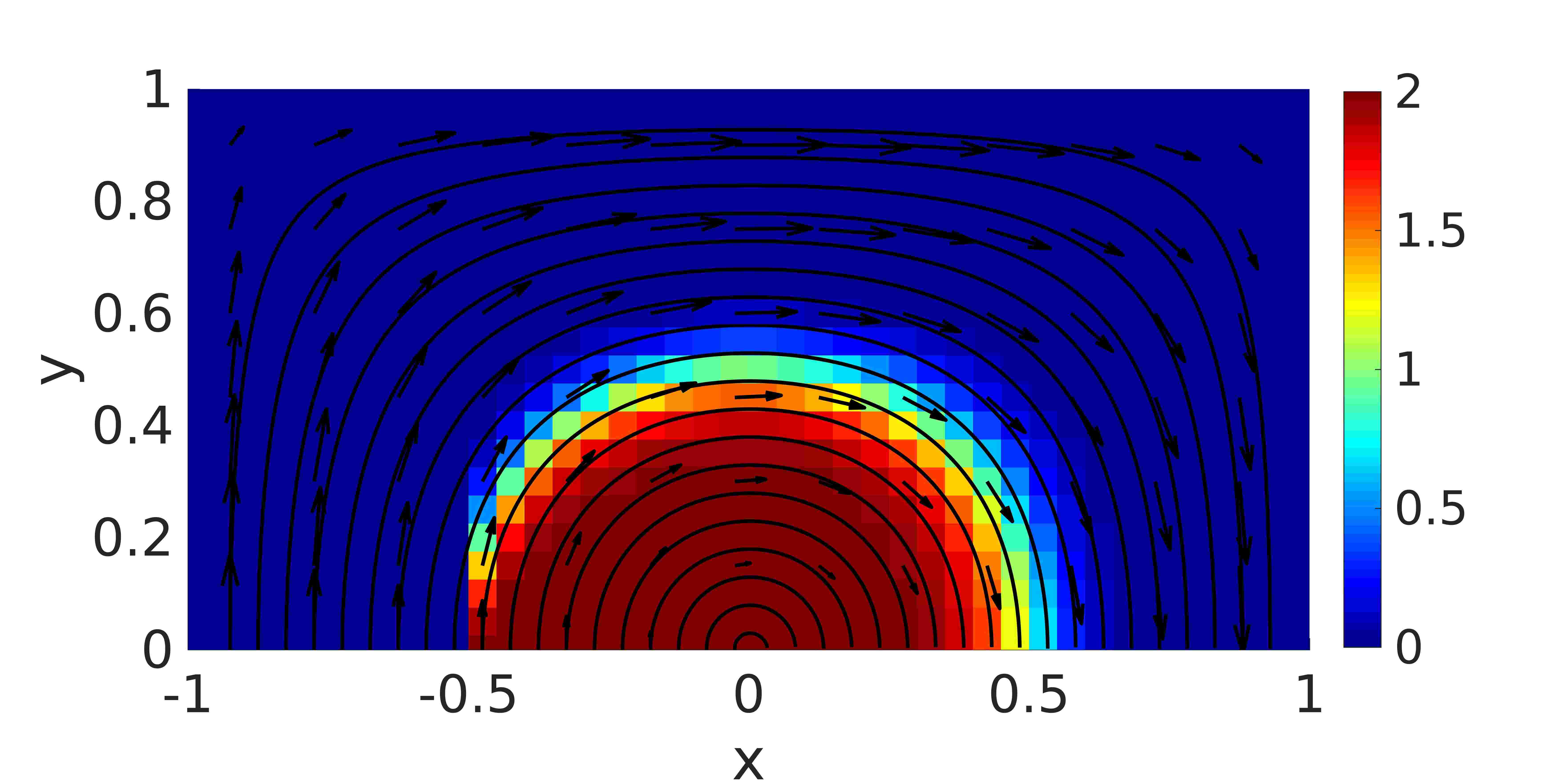}
			\end{subfigure}
			\begin{subfigure}[b]{1.0\textwidth}
				\rotatebox{90}{(m) QUICK(TVD)}
				\includegraphics[trim = 134 0 350 280,clip=true, keepaspectratio=true,width=0.9\textwidth]{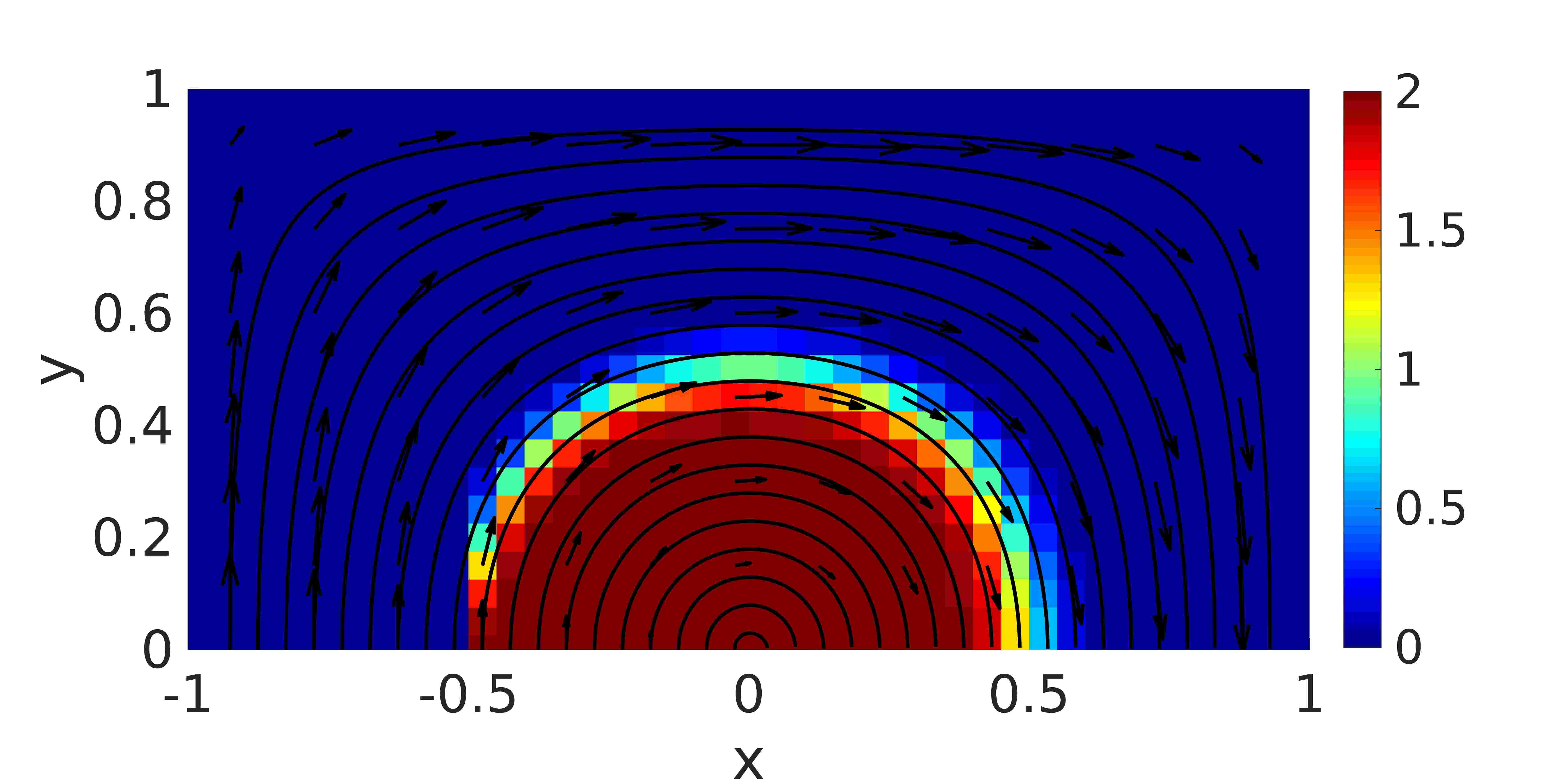}
			\end{subfigure}
			\begin{subfigure}[b]{1.0\textwidth}
				\rotatebox{90}{(n) SUPERBEE}
				\includegraphics[trim = 134 0 350 280,clip=true, keepaspectratio=true,width=0.9\textwidth]{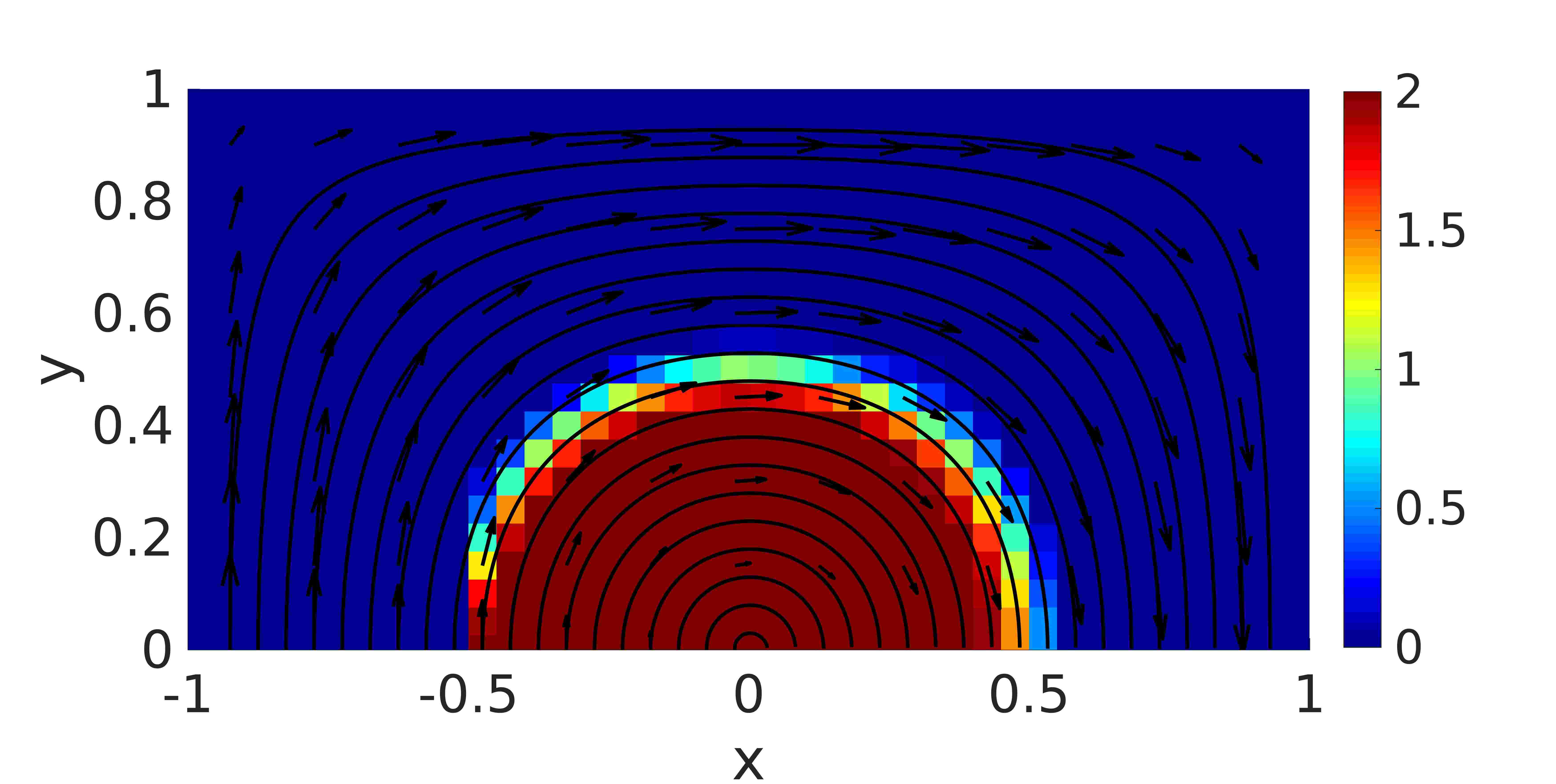}
			\end{subfigure}
			\ \\[3mm]
			\ \\
			\ \\
			\ \\
			\ \\
			\ \\
			\ \\
			\ \\
			\ \\
			\ \\
			\ \\
			\ \\
			\ \\
			\ \\
			\ \\
			\ \\
		\end{minipage}
	\end{minipage}
	\caption{Smith and Hutton problem contour plots of $\phi$ on a mesh with $40 \times 20$ nodes for $\alpha=1000$ and $\Gamma = 0$. Results were obtained using (a) - (e) DStreaM with non-limited coefficients (\ref{eq2.4}) with ranges from 1 to 5, respectively, upwind first-order scheme (f), (g) - (k) DStreaM with limited coefficients (\ref{eq2.5}) with ranges from 1 to 5, respectively, and TVD schemes with limiters Min-Mode (l), QUICK (m), and SUPERBEE (n).}
	\label{SH_contour_plot}
\end{figure}

\begin{figure}[htb!]
	\centering
	\includegraphics[height=0.3\textwidth]{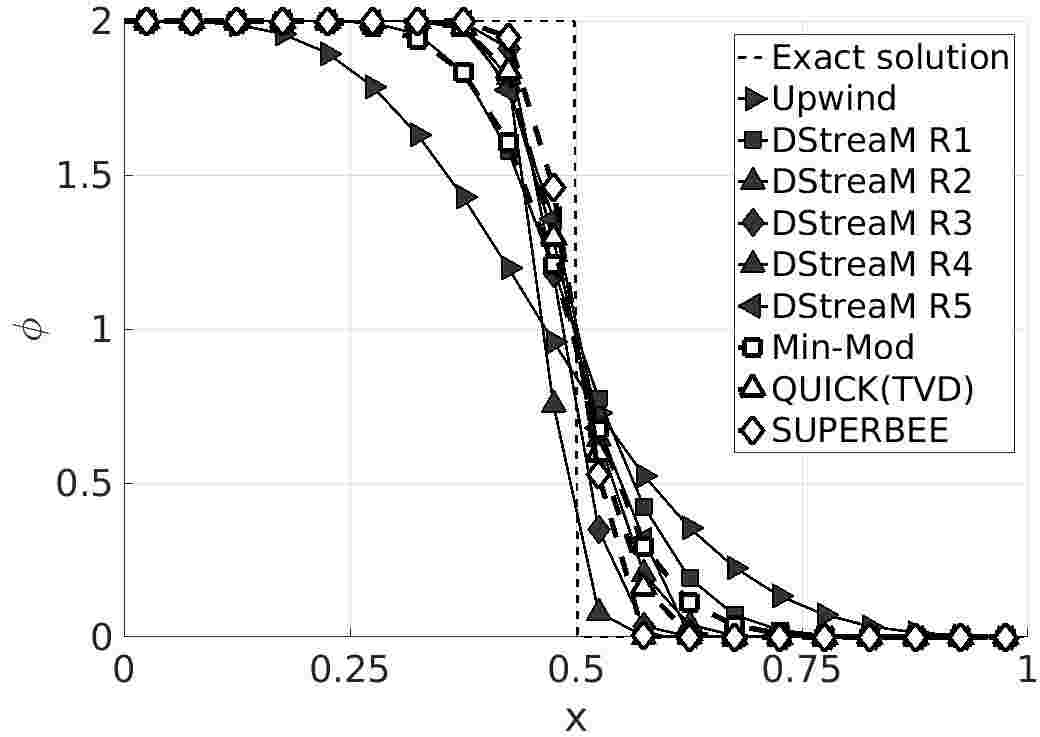}\\[-5mm]
	\caption{Obtained profiles at the outlet ($0 \le x \le 1$ and $y=0$) for Smith and Hutton problem, $\alpha=1000$, $\Gamma = 0$.\\[-15mm]}
	\label{Smith-Hutton_problem_alpha1000_Gamma0}
\end{figure}

\begin{figure}[htbp!]
	\centering
	\scriptsize
	\begin{minipage}{0.49\textwidth}
		\begin{subfigure}[b]{1.0\textwidth}
			\rotatebox{90}{(a) Upwind}
			\includegraphics[trim = 134 0 350 280,clip=true, keepaspectratio=true,width=0.9\textwidth]{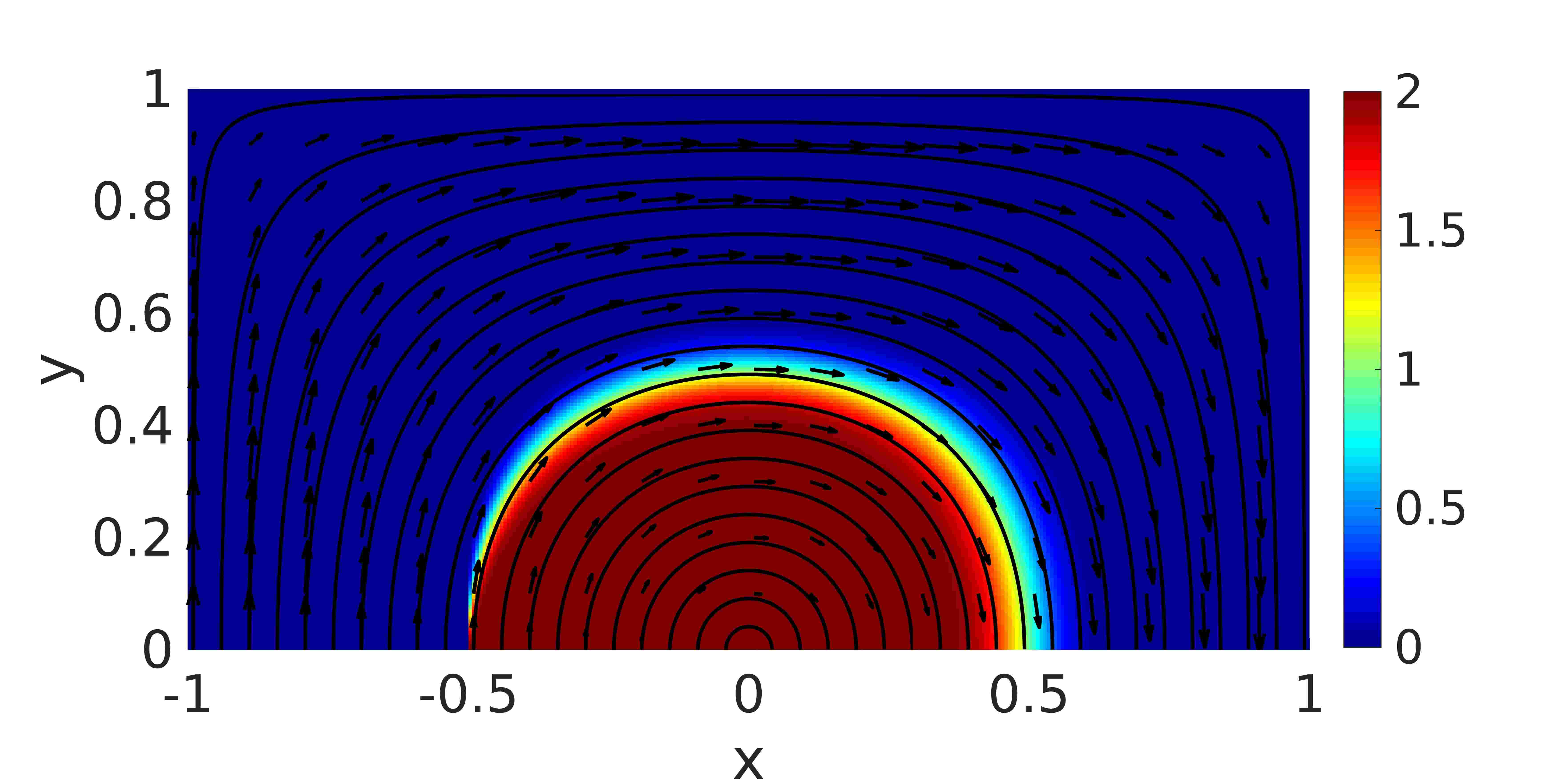}
			\label{SH_contour_plot_finer_mesh:UW}
		\end{subfigure}
		\begin{subfigure}[b]{1.0\textwidth}
			\rotatebox{90}{(b) DStreaM R1}
			\includegraphics[trim = 134 0 350 280,clip=true, keepaspectratio=true,width=0.9\textwidth]{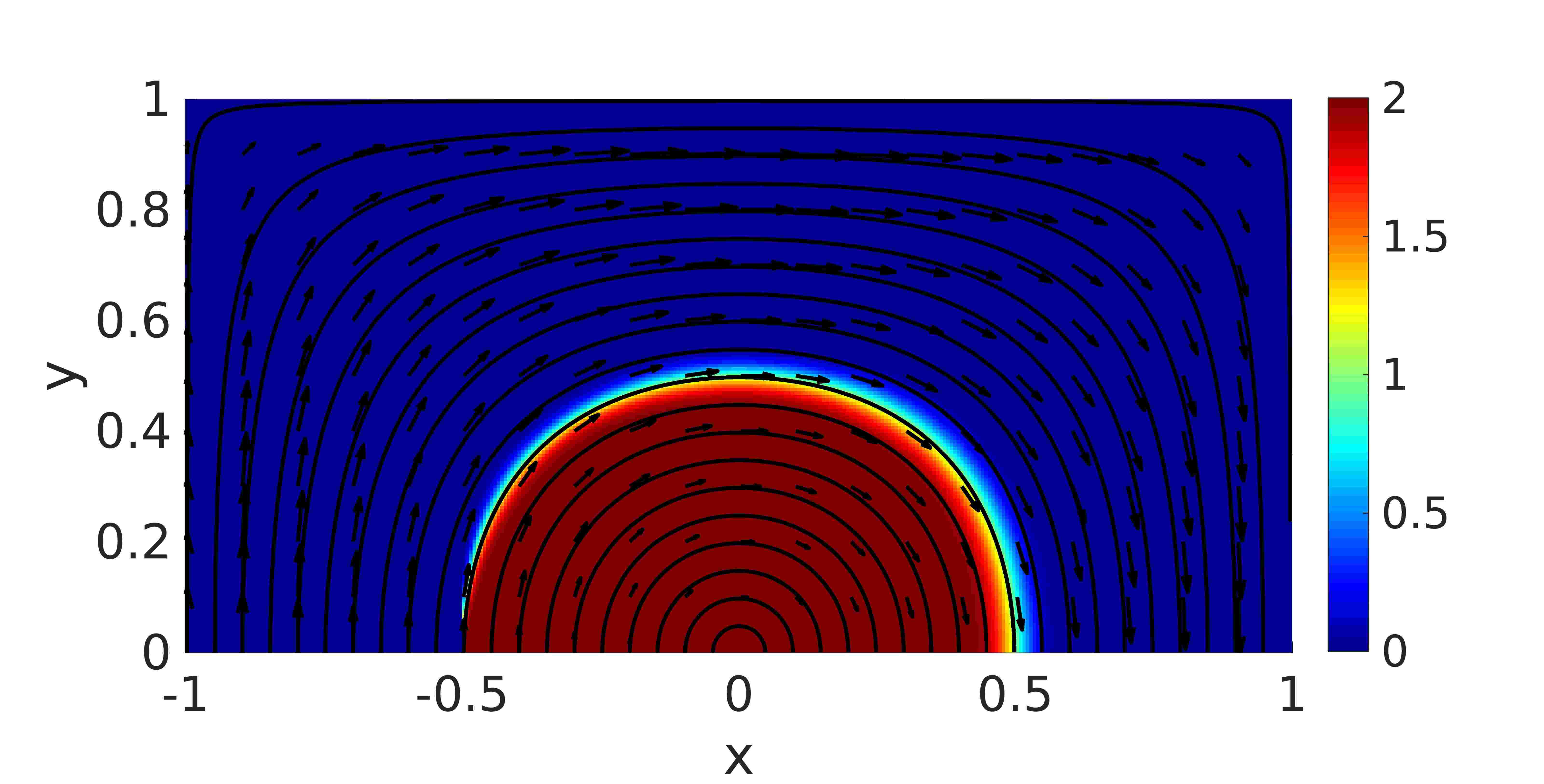}
			\label{SH_contour_plot_finer_mesh:UI1}
		\end{subfigure}
		\begin{subfigure}[b]{1.0\textwidth}
			\rotatebox{90}{(c) DStreaM R2}
			\includegraphics[trim = 134 0 350 280,clip=true, keepaspectratio=true,width=0.9\textwidth]{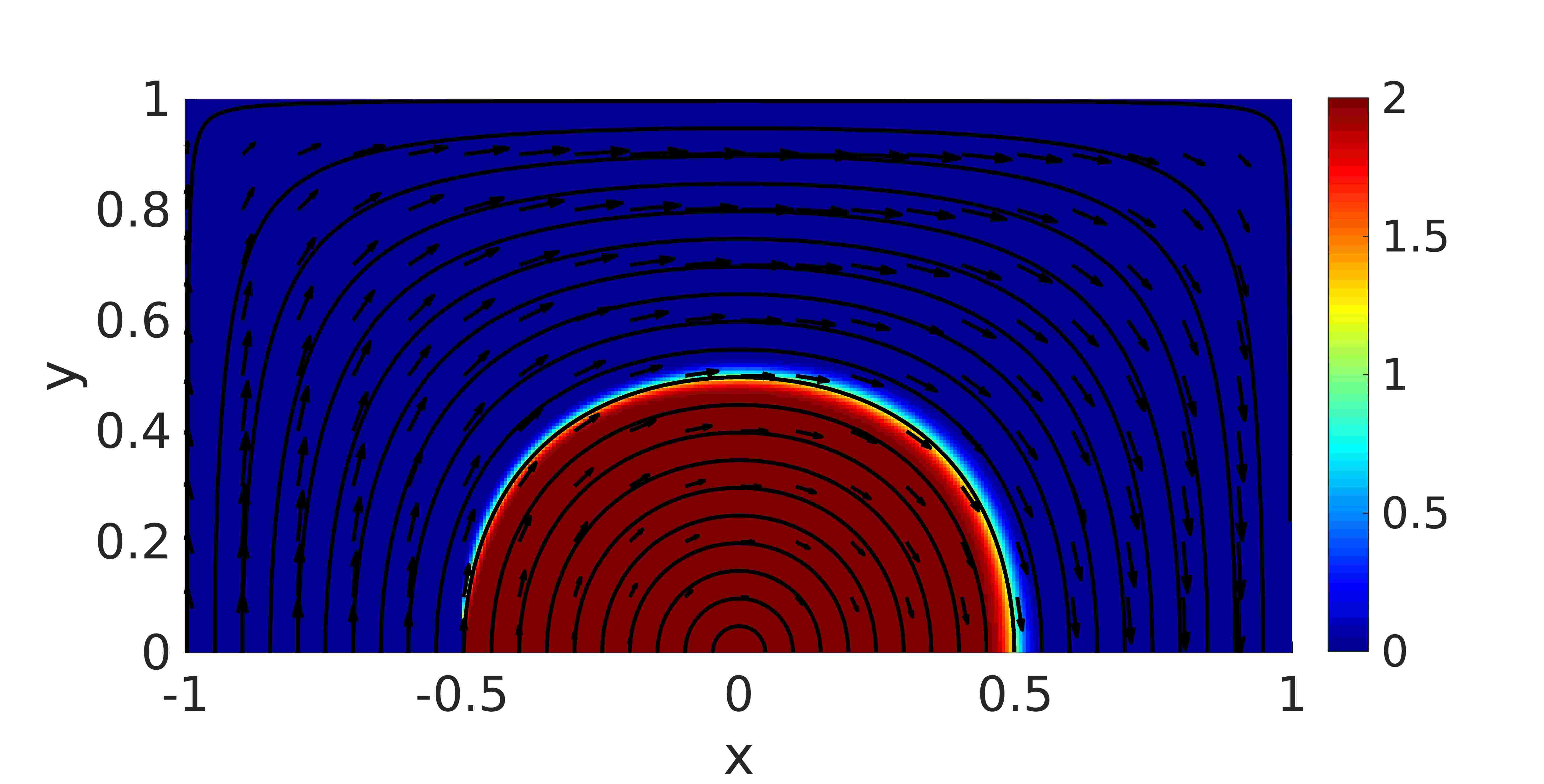}
			\label{SH_contour_plot_finer_mesh:UI2}
		\end{subfigure}
		\begin{subfigure}[b]{1.0\textwidth}
			\rotatebox{90}{(d) DStreaM R3}
			\includegraphics[trim = 134 0 350 280,clip=true, keepaspectratio=true,width=0.9\textwidth]{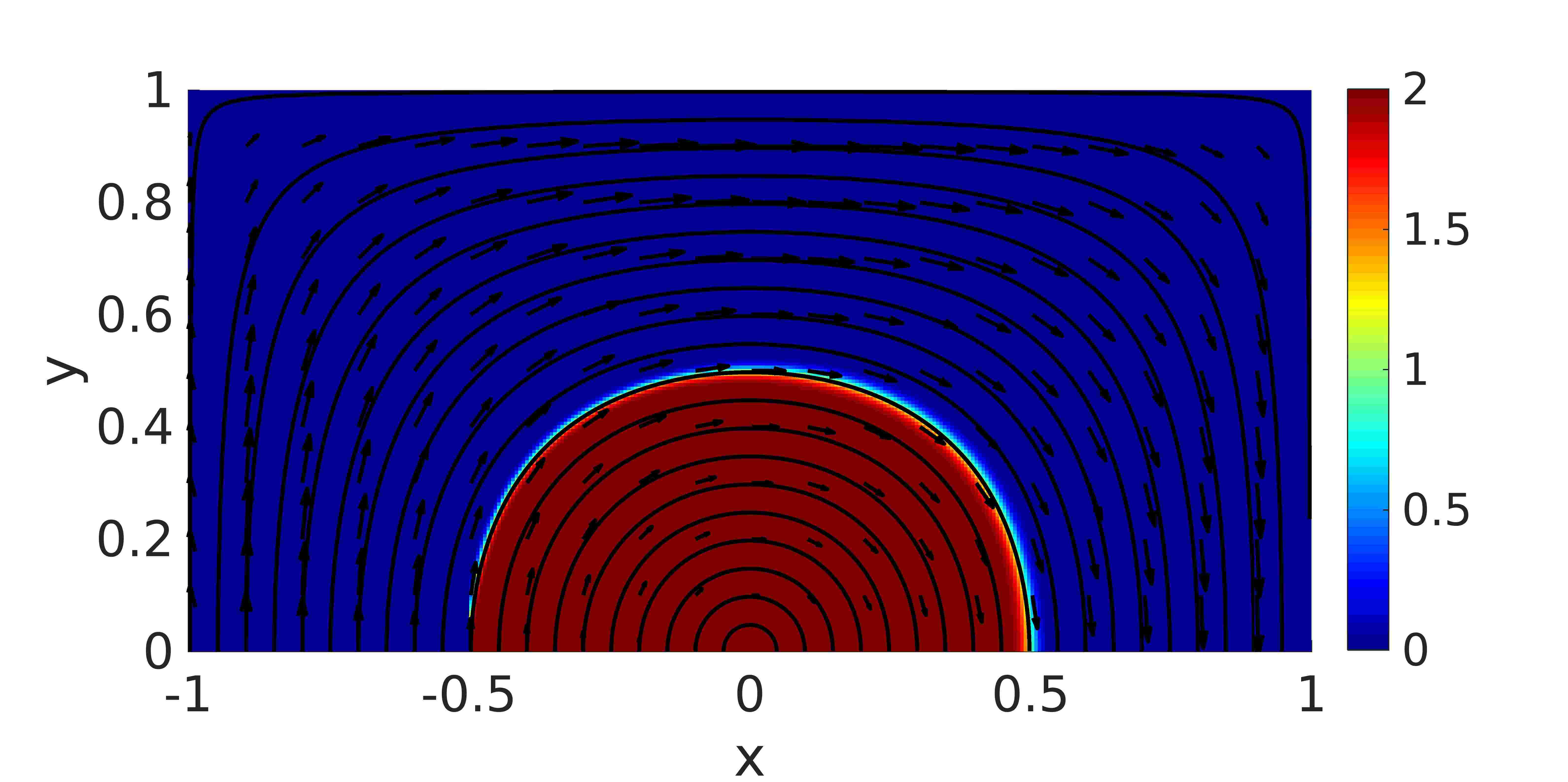}
			\label{SH_contour_plot_finer_mesh:UI3}
		\end{subfigure}
		\begin{subfigure}[b]{1.0\textwidth}
			\rotatebox{90}{(e) DStreaM R4}
			\includegraphics[trim = 134 0 350 280,clip=true, keepaspectratio=true,width=0.9\textwidth]{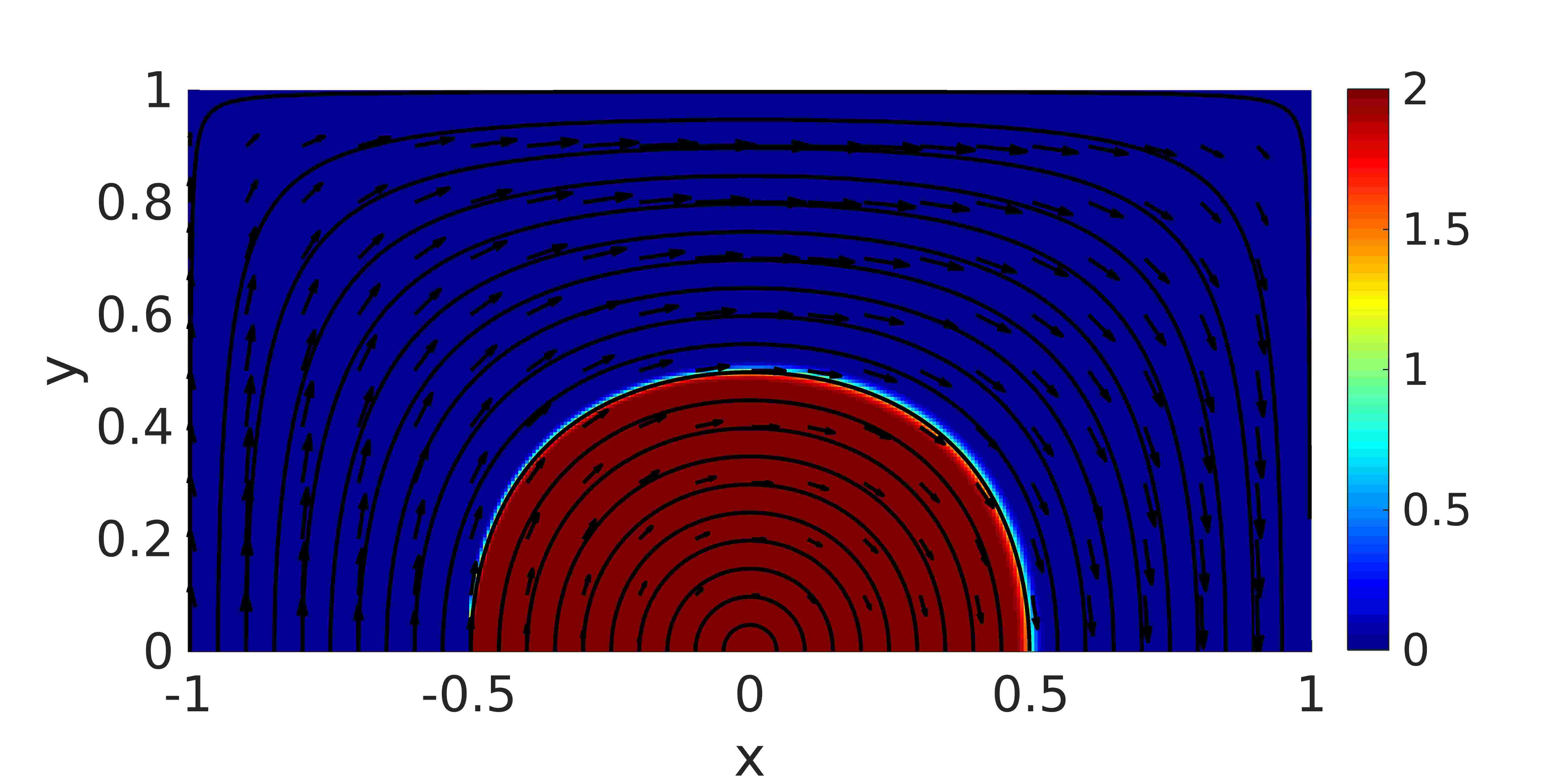}
			\label{SH_contour_plot_finer_mesh:UI4}
		\end{subfigure}
		\begin{subfigure}[b]{1.0\textwidth}
			\rotatebox{90}{(f) DStreaM R5}
			\includegraphics[trim = 134 0 350 280,clip=true, keepaspectratio=true,width=0.9\textwidth]{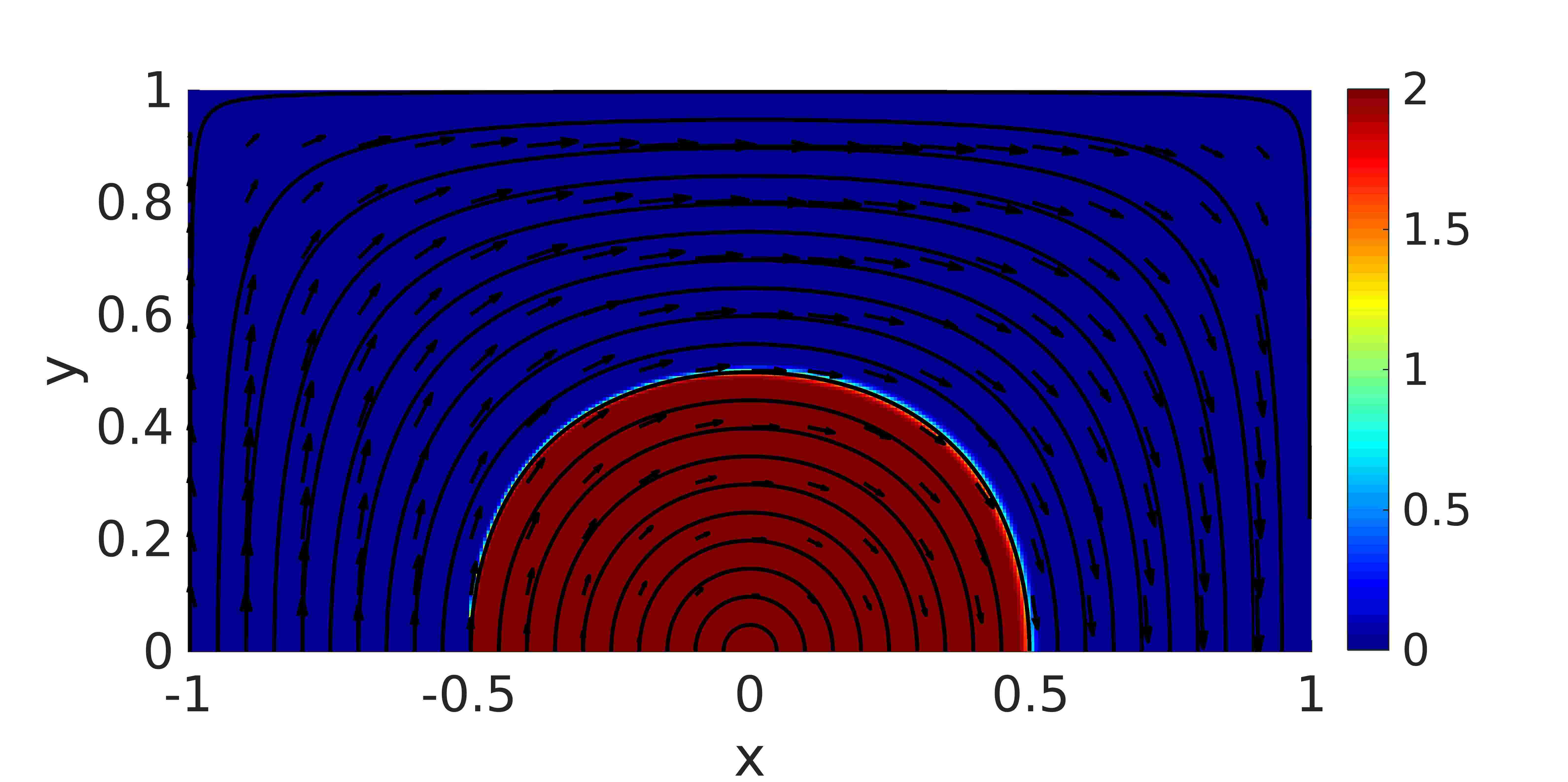}
			\label{SH_contour_plot_finer_mesh:UI5}
		\end{subfigure}
	\end{minipage}
	\begin{minipage}{0.49\textwidth}
		\ \\
		\ \\
		\ \\
		\ \\
		\ \\
		\ \\
		\ \\
		\ \\
		\ \\
		\ \\
		\ \\
		\ \\
		\ \\
		\ \\[-1mm]
		\begin{subfigure}[b]{1.0\textwidth}
			\rotatebox{90}{(g) Min-Mod}
			\includegraphics[trim = 134 0 350 280,clip=true, keepaspectratio=true,width=0.9\textwidth]{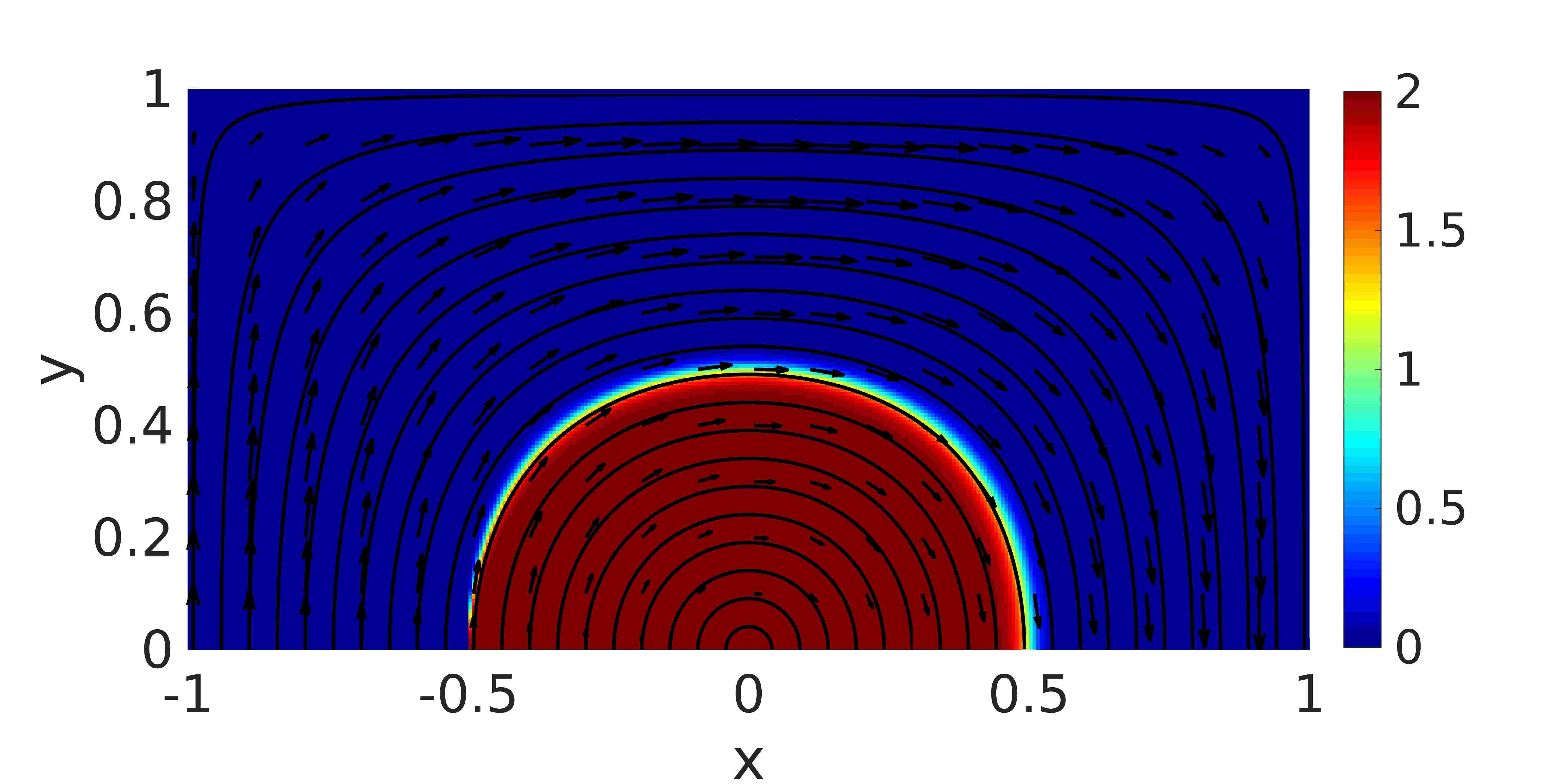}
			\label{SH_contour_plot_finer_mesh:Min-Mod}
		\end{subfigure}
		\begin{subfigure}[b]{1.0\textwidth}
			\rotatebox{90}{(h) QUICK(TVD)}
			\includegraphics[trim = 134 0 350 280,clip=true, keepaspectratio=true,width=0.9\textwidth]{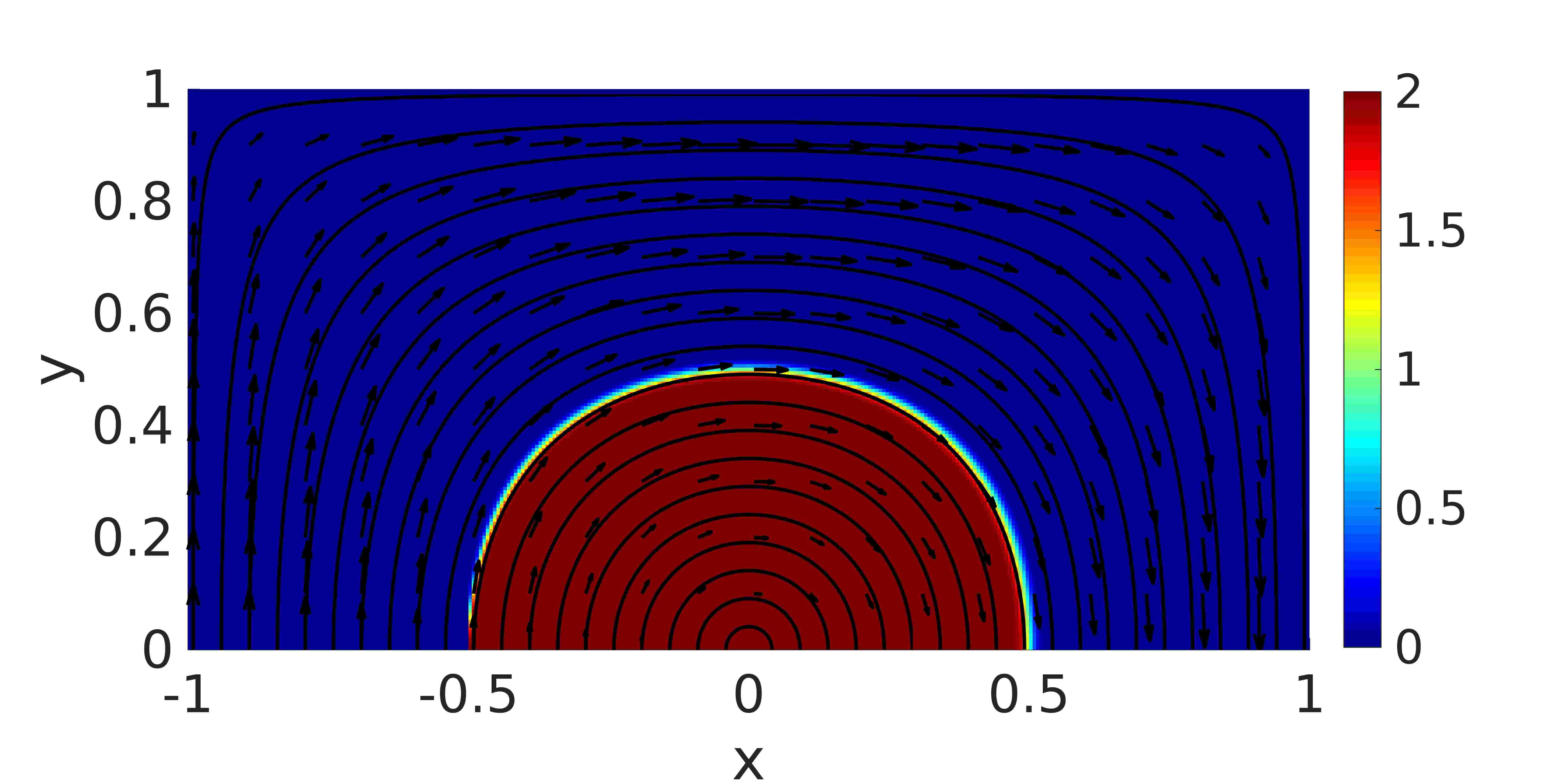}
			\label{SH_contour_plot_finer_mesh:QUICK(TVD)}
		\end{subfigure}
		\begin{subfigure}[b]{1.0\textwidth}
			\rotatebox{90}{(i) SUPERBEE}
			\includegraphics[trim = 134 0 350 280,clip=true, keepaspectratio=true,width=0.91\textwidth]{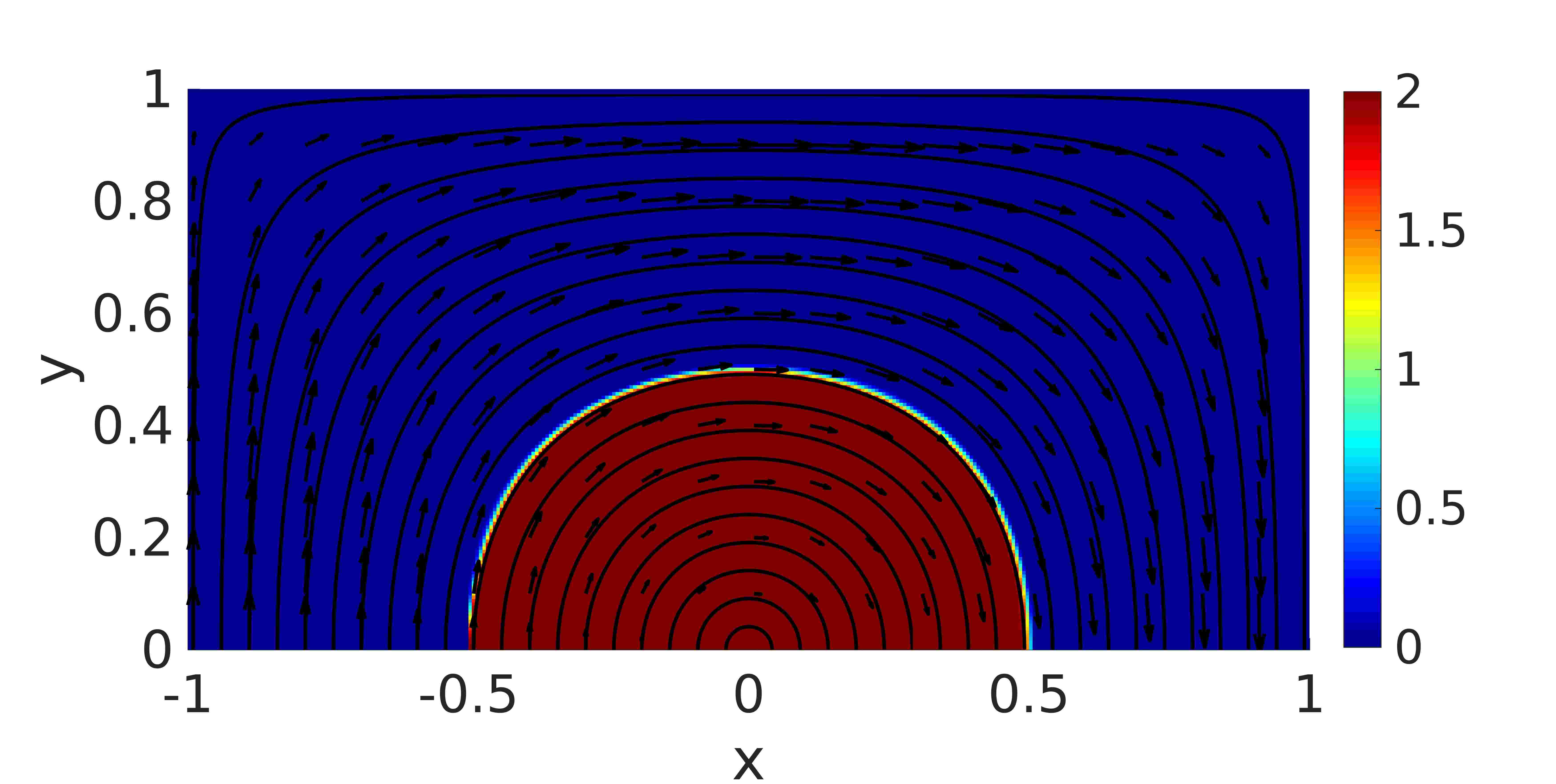}
			\label{SH_contour_plot_finer_mesh:SUPERBEE}
		\end{subfigure}
	\end{minipage}
	\caption{Smith and Hutton problem contour plots of $\phi$ on a mesh with $320 \times 160$ nodes for $\alpha=1000$ and $\Gamma = 0$. Results were obtained using upwind first-order scheme (a), (b) - (f) DStreaM with limited coefficients (\ref{eq2.5}) with ranges from 1 to 5, respectively, and TVD schemes with limiters Min-Mode (g), QUICK (h), and SUPERBEE (i).}
	\label{SH_contour_plot_finer_mesh}
\end{figure}

\begin{figure}[htb!]
\centering
\begin{subfigure}[b]{0.48\textwidth}
	\centering
	\includegraphics[width=\textwidth]{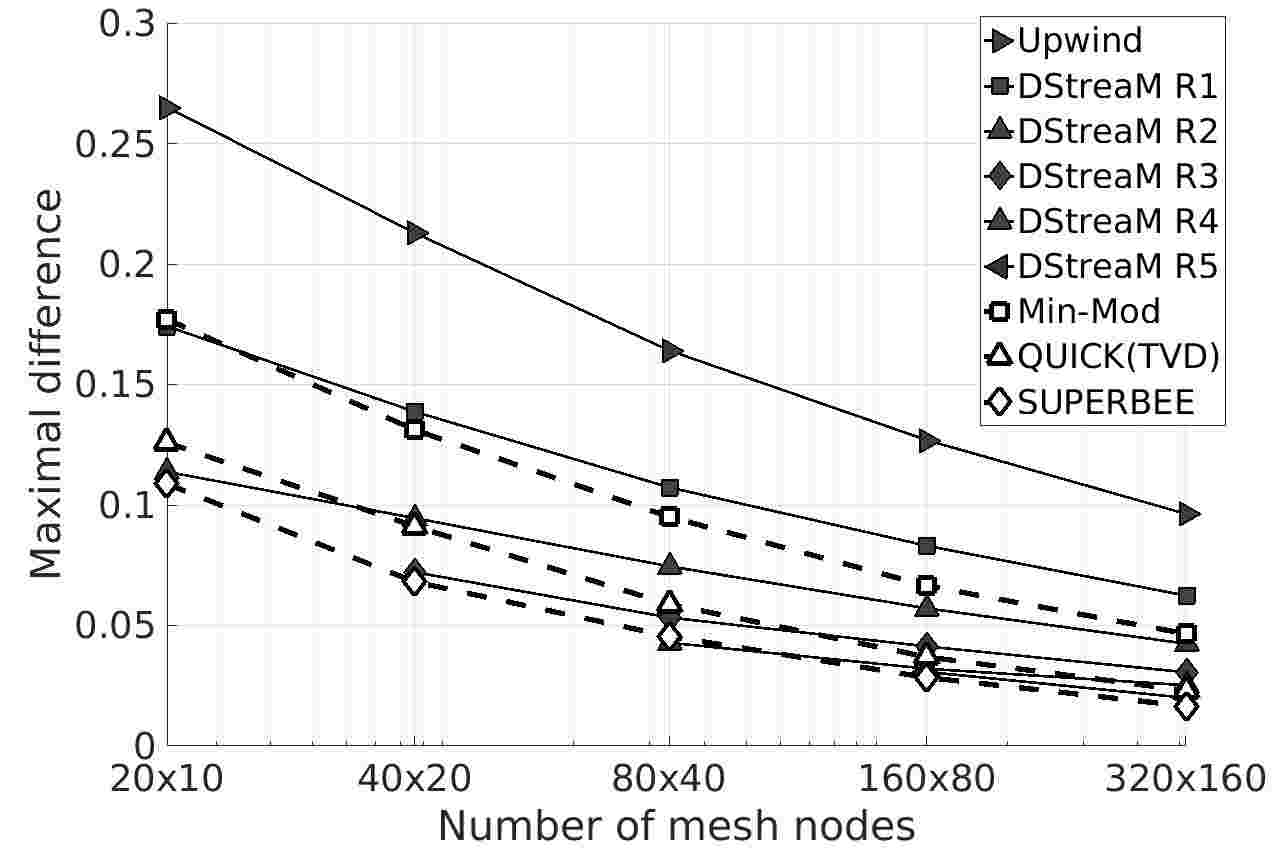}
	\caption{Mesh convergence.\\ \ \\[-5mm]}
	\label{Smith-Hutton_mesh_convergence:Max_Diff}
\end{subfigure}
\begin{subfigure}[b]{0.48\textwidth}
	\centering
	\includegraphics[width=\textwidth]{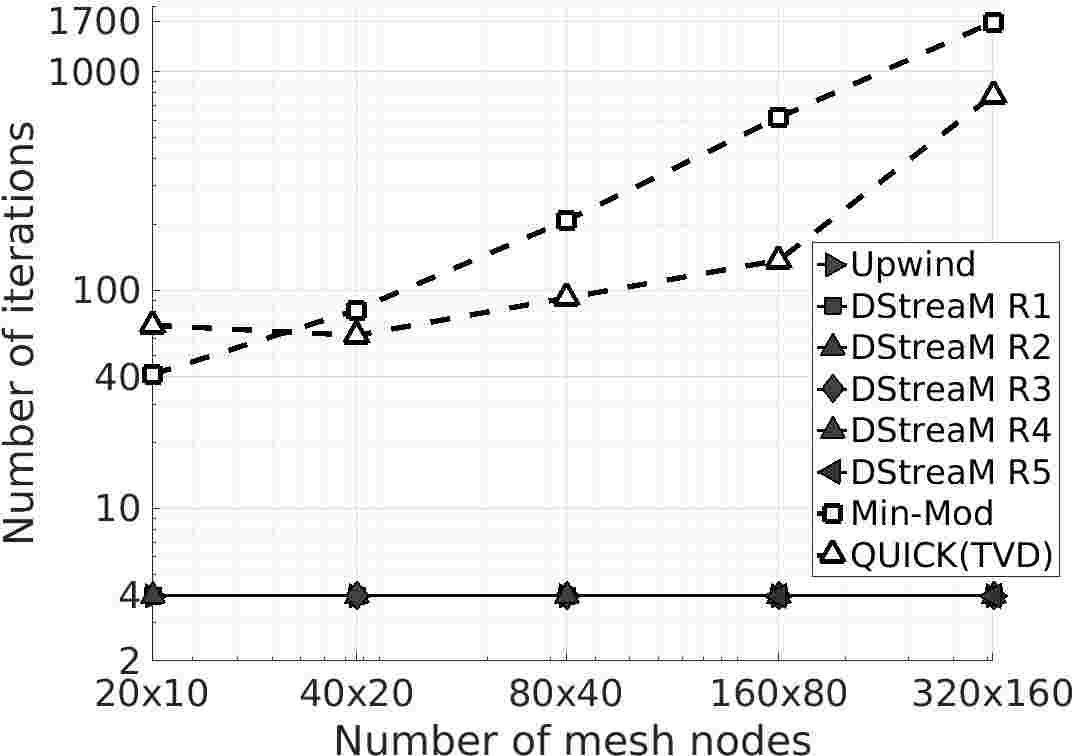}
	\caption{Number of iterations as function of number of mesh nodes.\\[-5mm]}
	\label{Smith-Hutton_mesh_convergence:N_Iter}
\end{subfigure}
	\caption{Smith and Hutton problem: (a) mesh convergence and (b) number of iterations as function of number of mesh nodes.\\[-8mm]}
	\label{Smith-Hutton_mesh_convergence}
\end{figure}

The Smith and Hutton problem \cite{Smith_and_Hutton_1982} will be investigated at $u=2y(1-x^2)$ and $v=-2x(1-y^2)$. Along the inlet $\phi$ is prescribed by the distribution of $\phi(-1 \leq x \leq 0, y=0) = 1+\text{tanh}[\alpha (2 x + 1)]$, see Fig. \ref{Smith-Hutton_problem}. Along the lines $x=-1$, $y=1$, and $x=1$, $\phi$ is prescribed as $1-\text{tanh}(\alpha)$ while along the outlet $(0 \leq x \leq 1, y=0)$ a zero gradient is specified for $\phi$. The considered case is for $\alpha=1000$ and $\Gamma=0$. Detailed comparisons are presented for mesh with $40 \times 20$ nodes. After that, approximation schemes mesh convergence is presented.\\\indent
As was considered earlier DStreaM with limited coefficients (\ref{eq2.5}) obtains results without under- or over-shoots (Fig. \ref{SH_contour_plot} (g)-(k)) while DStreaM with non-limited coefficients (\ref{eq2.4}) obtains results with under/over-shoots when the range is equal to or greater than 4, (Fig. \ref{SH_contour_plot} (a)-(e)). The reason is that local streamline through $\phi_{i,j}$ is curved and deviate from discrete streamline when the range is larger because when the range is larger discrete streamline approaches a straight line. In general, a curved line can be approximated within sufficient accuracy using small straight lines successfully. When we use the smaller spatial step as a case with $320 \times 160$ nodes, obtained results are without under/over-shoots for greater ranges than 4, see Fig. \ref{SH_contour_plot_finer_mesh} (a)-(f). The reason is that the range 5 for mesh with $320 \times 160$ nodes corresponds to distance $5 \Delta_{320 \times 160} = 5 * 1/160 = 0.03125$ ($\Delta_{320 \times 160}$ is the spatial step for mesh with $320 \times 160$ nodes) that is smaller than distance for mesh with $40 \times 20$ nodes that is $5 \Delta_{40 \times 20} = 5 * 1/20 = 0.25$ ($\Delta_{40 \times 20}$ is the spatial step for mesh with $40 \times 20$ nodes). Finer mesh makes discrete streamline length shorter and increases their number. That lead to a better approximation of curved streamline (non-uniform velocity field). That property of DStreaM makes it suitable for applications on unstructured meshes where smaller spatial steps are used in places of high gradients.\\\indent
The comparison shows that DStreaM spatial accuracy is second order as TVD schemes and requires 4 iterations as upwind scheme. Fig. \ref{SH_contour_plot} shows comparisons of obtained fields by DStreaM, upwind and TVD schemes with limiters Min-Mod, QUICK, and SUPERBEE. DStreaM results with ranges from 1 to 3 (Fig. \ref{SH_contour_plot} (g),(h), and (i)) correspond to results of TVD schemes with limiters Min-Mod, QUICK, and SUPERBEE (Fig. \ref{SH_contour_plot} (l), (m), and (n)), respectively. Profiles of obtained solutions along the line $(0 \leq x \leq 1, y=0)$ are plotted in Fig. \ref{Smith-Hutton_problem_alpha1000_Gamma0}. The comparison shows second order accuracy of the proposed approach DStreaM. Upwind and DStreaM obtain a final solution for 4 iterations while Min-Mod and QUICK need 81 and 62 iterations, see Table \ref{Smith-Hutton_problem_Performance_40x20_10_8}. TVD schemes require approximately 20 and 15 times more iterations to obtain the same solution as DStreaM. SUPERBEE scheme cannot obtain a solution within considered residual despite small under-relaxation coefficient and many iterations while DStreaM R3 can obtain an as sharp solution as SUPERBEE for four iterations.
\newgeometry{left=2cm,right=1.5cm,top=2cm,bottom=2cm}
\captionsetup[subfigure]{position=bottom, labelfont=bf,textfont=normalfont,singlelinecheck=off,justification=raggedright}
\restoregeometry
\captionsetup[subfigure]{position=bottom, labelfont=bf,textfont=normalfont,singlelinecheck=off,justification=raggedright}
\begin{table}[h]
\caption{Under-relaxation coefficient, number of iterations, and obtained a maximum residual of considered schemes for calculation of Smith and Hutton problem with convergence criteria $\epsilon = 10^{-8}$ on a mesh with $40 \times 20$ nodes .}
\centering
\begin{tabular}{lrrrr}
\toprule
Approximation scheme & $\alpha$ & Iterations & max residual \\
\midrule
Upwind            & 1.0  &     4 & 0 \\
DStreaM R1        & 1.0  &     4 & 0 \\
DStreaM R2        & 1.0  &     4 & 0 \\
DStreaM R3        & 1.0  &     4 & 0 \\
DStreaM R4        & 1.0  &     4 & 0 \\
DStreaM R5        & 1.0  &     4 & 0 \\
Min-Mod           & 0.95 &    81 & $5.5 \times 10^{-9}$ \\
QUICK(TVD)        & 0.85 &    62 & $9.0 \times 10^{-9}$ \\
SUPERBEE          & 0.1  & 50000 & $10^{-3}$ \\
\bottomrule
\end{tabular}
\label{Smith-Hutton_problem_Performance_40x20_10_8}
\end{table}
Mesh convergence study shows DStreaM advantages. Fig. \ref{Smith-Hutton_mesh_convergence} (a) shows a maximal difference between the numerical and exact solution of obtained fields as a function of the number of mesh nodes. Presented solutions obtained by DStreaM are without unphysical oscillations (when limited coefficients are not needed). All schemes reduced obtained maximal difference for finer mesh as TVD schemes demonstrate better mesh convergence compared to DStreaM with a fixed range. Figures \ref{SH_contour_plot} and \ref{SH_contour_plot_finer_mesh} also show DStreaM accuracy "shifting". DStreaM can use larger ranges on a finer mesh to reduce maximal difference and to keep its spacial accuracy closer to SUPERBEE. Fig. \ref{Smith-Hutton_mesh_convergence} (b) presents a number of iterations needed to obtain a solution with convergence criteria $\epsilon=10^{-8}$. TVD schemes required under-relaxation coefficient to converge iterative solver. As a number of iterations vary according to under-relaxation coefficients, results presented in Fig. \ref{Smith-Hutton_mesh_convergence} (b) for TVD schemes is minimally obtained for a range of under-relaxation coefficients. SUPERBEE obtains the solution within required convergence criteria only for mesh with $20 \times 10$ nodes for 68 iterations while for the others meshes it reaches maximum residual $\sim10^{-3}$, and the number of iterations for SUPERBEE are not presented. On the other hand, upwind and DStreaM do not require under-relaxation coefficients to ensure convergence of iterative process. TVD schemes require more iterations on finer mesh and number of iterations rapidly increases from 40 to 1700. Contrary, upwind and DStreaM need 4 iterations that are independent from mesh refinement. The results show that TVD schemes demonstrate better mesh convergent than DStreaM with a fixed range, but DStreaM can increase ranges for finer mesh and to obtain results with corresponding spatial accuracy to SUPERBEE, see Fig. \ref{SH_contour_plot_finer_mesh}. After all, DStreaM needs 4 iterations to obtain a final solution that is from 10 to 425 times fewer iterations compared to TVD schemes, and it can obtain an as sharp solution as SUPERBEE within a given residual when range varies.
\FloatBarrier
\section{Conclusions}
In this paper was proposed Discrete Stream(line) Method (DStreaM) for convective terms approximation. The main idea is the approximation approach to corresponds to a physical phenomenon described by approximated terms. In considered case, the physical phenomenon is pure convection that is described by convective terms. Convective terms describe transporting a property along the streamline, and the main characteristics are unidirectional information propagation without exchange information with transported properties over neighbor streamlines. DStreaM approach maps streamlines to a mesh. As a result, streamlines are represented as Discrete Streamlines. Discrete Streamlines represent unidirectional information propagation in mesh orientated according to velocity field. Four standard advection (pure convection) test problems were presented in the paper: advection of a step profile, advection of a double-step profile, advection of a sinusoidal profile, and Smith and Hutton problem. DStreaM's results were compared with results obtained by upwind first-order scheme and TVD schemes with limiters Min-Mod, QUICK, and SUPERBEE. TVD schemes use under-relaxation coefficients to ensure convergence and require from 15 to 93.5 times more iterations to obtain its final solution than upwind and DStreaM. For Smith and Hutton problem SUPERBEE cannot obtain a solution within given convergence criteria, DStreaM can obtain an as sharp solution as SUPERBEE within 4 iterations. DStreaM combines the positives of upwind and TVD schemes. It obtains a solution within 2 or 4 iterations without under-relaxation coefficients as upwind, and it is second order scheme as TVD schemes.\\\indent
DStreaM approach looks promising for calculation of convective-dominated problems because:
\begin{itemize}
 \item Naturally approximates first derivatives in partial differential equations;
 \item All independent variables can be defined in the same nodes;
 \item It can be straightforward and easy applied on arbitrary meshes as a meshfree method or unstructured meshes with real-time mesh refinement;
\end{itemize}

\section*{Acknowledgments}
We would like to acknowledge the financial support provided by the Bulgarian NSF under Grant DN-02/7-2016.

\bibliographystyle{siamplain}
\bibliography{bibliography_list_1}
\end{document}